\documentclass[a4paper, english, 11pt, reqno, dvipsnames]{amsart}

\usepackage{amssymb}

\usepackage[english]{babel}
\usepackage[utf8]{inputenc}
\usepackage{tikz,tikzscale}
\tikzset{>=latex}
\usetikzlibrary{patterns}
\usepackage{enumerate,enumitem}
\usepackage{booktabs,multirow}
\usepackage{amsfonts,mathtools}
\usepackage{dsfont}
\usepackage[square,numbers]{natbib}

\usepackage{pgfplots}
\pgfplotsset{compat=1.13}

\usepackage{hyperref}
\usepackage[noabbrev, capitalise, nameinlink]{cleveref}
\crefname{equation}{}{}
\crefname{assumption}{Assumption}{Assumptions}
\crefformat{equation}{\textup{#2(#1)#3}}

\newtheorem{theorem}{Theorem}[section]

\newtheorem{lemma}[theorem]{Lemma}

\theoremstyle{definition}

\theoremstyle{remark}
\newtheorem{remark}[theorem]{Remark}
\numberwithin{theorem}{section}
\numberwithin{equation}{section}
\numberwithin{table}{section}
\numberwithin{figure}{section}

\renewcommand*{\vec}{\boldsymbol}

\newcommand*{\setEdge}{\ensuremath{\mathcal E}}
\newcommand*{\setNode}{\ensuremath{\mathcal N}}
\newcommand*{\edge}{\ensuremath{\mathfrak e}}
\newcommand*{\node}{\ensuremath{\mathfrak n}}

\let\rho\varrho

\newcommand{\setNodeDir}{\ensuremath{\setNode_\textup D}}

\newcommand{\jump}[1]{{[\![ #1 ]\!]}}

\newcommand{\IR}{\ensuremath{\mathbb R}}

\DeclareMathAlphabet{\mathbfsf}{\encodingdefault}{\sfdefault}{bx}{n}

\newcommand{\ds}{\ensuremath{\, \textup d \sigma}}

\newcommand{\re}{\ensuremath{\vec r_\edge}}
\newcommand{\rn}{\ensuremath{\vec r_\node}}
\newcommand{\ue}{\ensuremath{\vec u_\edge}}
\newcommand{\un}{\ensuremath{\vec u_\node}}
\newcommand{\und}{\ensuremath{\vec u_\node^\textup{D}}}
\newcommand{\rnd}{\ensuremath{\vec r_\node^\textup{D}}}
\newcommand{\redisc}{\ensuremath{\bar{\vec r}_\edge}}
\newcommand{\rndisc}{\ensuremath{\bar{\vec r}_\node}}
\newcommand{\uedisc}{\ensuremath{\bar{\vec u}_\edge}}
\newcommand{\undisc}{\ensuremath{\bar{\vec u}_\node}}
\newcommand{\ndisc}{\ensuremath{\bar{\vec n}_\edge}}
\newcommand{\mdisc}{\ensuremath{\bar{\vec m}_\edge}}
\newcommand{\Ndisc}{\ensuremath{\bar{\vec N}_\edge}}
\newcommand{\Mdisc}{\ensuremath{\bar{\vec M}_\edge}}
\newcommand{\Udisc}{\ensuremath{\bar{\vec U}_\edge}}
\newcommand{\Rdisc}{\ensuremath{\bar{\vec R}_\edge}}
\newcommand{\lambdadisc}{\ensuremath{\bar{\vec \lambda}}}
\newcommand{\phidisc}{\ensuremath{\bar{\vec \phi}}}
\newcommand{\mudisc}{\ensuremath{\bar{\vec \mu}}}
\newcommand{\psidisc}{\ensuremath{\bar{\vec \psi}}}

\newcommand{\with}{\,:\,}

\newcommand{\uhat}{\ensuremath{\hat{\vec u}}}
\newcommand{\rhat}{\ensuremath{\hat{\vec r}}}
\newcommand{\nhat}{\ensuremath{\hat{\vec n}}}
\newcommand{\mhat}{\ensuremath{\hat{\vec m}}}
\newcommand{\ihat}{\ensuremath{\hat{\vec i}}}
\newcommand{\lamo}{\ensuremath{\hat{\vec \lambda}_1}}
\newcommand{\lamt}{\ensuremath{\hat{\vec \lambda}_2}}
\newcommand{\phio}{\ensuremath{\hat{\vec \phi}_1}}
\newcommand{\phit}{\ensuremath{\hat{\vec \phi}_2}}
\newcommand{\lamdiff}{\ensuremath{\hat{\vec \lambda}_\Delta}}
\newcommand{\phidiff}{\ensuremath{\hat{\vec \phi}_\Delta}}
\newcommand{\phisum}{\ensuremath{\hat{\vec \phi}_\Sigma}}
\newcommand{\edgelength}{\ensuremath{h_\mathfrak{e}}}
\newcommand{\maxedgelength}{\ensuremath{h}}

\makeatletter
\newcommand\footnoteref[1]{\protected@xdef\@thefnmark{\ref{#1}}\@footnotemark}
\makeatother

\begin{document}

\title[HDG for Timoshenko beam network models ]{Arbitrary order approximations at constant cost for Timoshenko beam network models } 

\author[M.~Hauck]{Moritz Hauck}
\author[A.~M{\aa}lqvist]{Axel M{\aa}lqvist}
\address{Department of Mathematical Sciences, University of Gothenburg and Chalmers University of Technology, 41296 Göteborg, Sweden}
\email{hauck@chalmers.se,\;axel@chalmers.se}

\author[A.~Rupp]{Andreas Rupp}
\address{School of Engineering Science, Lappeenranta--Lahti University of Technology, P.O. Box 20, 53851 Lappeenranta, Finland}
\email{andreas@rupp.ink}

\thanks{M.~Hauck has been supported by the Knut and Alice Wallenberg foundation postdoctoral program in mathematics for researchers from outside Sweden, grant number KAW 2022.0260 and A.~Målqvist by the Swedish research council, grant number 2023-03258\_VR. A.\ Rupp has been supported by the Academy of Finland's grant number 350101 \emph{Mathematical models and numerical methods for water management in soils}, grant number 354489 \emph{Uncertainty quantification for PDEs on hypergraphs}, grant number 359633 \emph{Localized orthogonal decomposition for high-order, hybrid finite elements}, Business Finland's project number 539/31/2023 \emph{3D-Cure: 3D printing for personalized medicine and customized drug delivery}, and the Finnish \emph{Flagship of advanced mathematics for sensing, imaging and modeling}, decision number 358944.}


\begin{abstract}
 This paper considers the numerical solution of Timoshenko beam network models, comprised of Timoshenko beam equations on each edge of the network, which are coupled at the nodes of the network using rigid joint conditions. Through hybridization, we can equivalently reformulate the problem as a symmetric positive definite system of linear equations posed on the network nodes. This is possible since the nodes, where the beam equations are coupled, are zero-dimensional objects. To discretize the beam network model, we propose a hybridizable discontinuous Galerkin method that can achieve arbitrary orders of convergence under mesh refinement without increasing the size of the global system matrix. As a preconditioner for the typically very poorly conditioned global system matrix, we employ a two-level overlapping additive Schwarz method. We prove uniform convergence of the corresponding preconditioned conjugate gradient method under appropriate connectivity assumptions on the network. Numerical experiments support the theoretical findings of this work. 
\end{abstract}

\date{\today}
\maketitle

\section{Introduction}
%
Many applications in science and engineering involve geometrically complex structures composed of slender, effectively one-dimensional objects. Examples include blood vessels~\cite{FKOWW22}, porous materials~\cite{CEPT12}, or fiber-based materials~\cite{KMM20}. For such problems, resolving all microscopic details in a three-dimensional computer simulation can be computationally very demanding. Therefore, in many cases, it seems appropriate to describe the geometry by a spatial network represented by a graph $\mathcal{G}=(\mathcal{N},\mathcal{E})$ of nodes and edges which is embedded into a bounded domain $\Omega \subset \mathbb R^3$. The resulting spatial network model then involves one-dimensional differential equations on each edge, coupled by algebraic constraints at the nodes of the graph. In the following, we consider the elastic deformation of fiber-based materials, such as paper or cardboard, as a model problem. The spatial network underlying this problem is constructed as follows: Nodes are placed at the intersections of fibers, and an edge connects two nodes if a fiber is connecting them. Depending on the intersection area of the two fibers, additional nodes and edges may be added to strengthen the connection between the fibers. To obtain a spatial network model, we equip each edge with a Timoshenko beam model, cf.~\cite{Timoshenko1921,Carrera2011}, describing the elastic deformation of the corresponding fiber. A well-posed problem is then obtained by enforcing the continuity of displacements and rotations and the balance of forces and moments at the nodes of the spatial network, cf.~\cite{Lagnese1994}. For more details on this model, particularly the construction of the graph representing the geometric structure of cardboard, we refer to~\cite{KMM20,Grtz2022,Grtz2024}. 

For practically simulating spatial network models, one needs to perform a discretization of the one-dimensional Timoshenko beam equations posed on the edges. A popular approach for this are beam finite elements, which differ mainly in the number and placement of the degrees of freedom; see, e.g., \cite{Kapur1966,Davis1972,Thomas1975,Lees1982}. When discretizing the Timoshenko beam equations using beam elements, one often observes a shear-locking effect, which leads to an underestimation of the displacements. This effect is caused by underresolution and occurs mainly for beam elements with few degrees of freedom; see, e.g., the theoretical study in~\cite{Mukherjee2001}. Shear locking can be avoided by increasing the local degrees of freedom by subdividing the edges and using multiple beam elements, or by using higher-order beam elements. Note that shear locking can also be reduced by using numerical tricks such as under-integration, cf.~\cite{Prathap1982,Yokoyama1994}. An alternative to beam elements are so-called analytical ansatz functions, which are derived by analytically solving the Timoshenko beam equations on each edge; see, e.g., \cite{Reddy1997,Jeleni2009,KufnerLSSS18}. However, to find an explicit expression for the analytical solution, it is usually necessary to impose additional modeling assumptions, such as that the material coefficients are constant or that the distributed loads and moments are homogeneous or constant. 

In this paper, we apply a hybridizable discontinuous Galerkin (HDG) method to discretize the spatial network model under consideration. The use of such a discretization is motivated by the study in~\cite{RuppGK22}, where the authors rigorously define and discretize diffusion-type problems on networks of hypersurfaces. There, the authors claim that a hybrid formulation is natural for partial differential equations on hypergraphs and, by extension, they should be discretized using HDG methods. In the special case of graphs, the nodes at which the one-dimensional differential equations are coupled are zero-dimensional objects, whence the spatial network problem can be equivalently reformulated as a symmetric positive definite system of linear equations posed on the nodes of the network. For such problems, an HDG discretization can achieve arbitrary convergence orders without increasing the number of globally coupled degrees of freedom.  Following the paradigm of HDG methods that local solves are essentially for free since the corresponding problems are small and they can be solved in parallel, cf.~\cite{CockburnGL09}, the proposed method can achieve arbitrary orders of convergence at (almost) constant computational cost. We perform an a priori error analysis of the HDG method, where we prove optimal convergence orders under mesh refinement. Note that the shear-locking effect of classical beam elements can be easily avoided by using sufficiently high polynomial degrees.  

Due to the complex geometry of the spatial network and possibly highly varying material coefficients, the linear system of equations obtained by the HDG method is typically very poorly conditioned. Numerical experiments demonstrate that standard black-box preconditioners, such as many algebraic multigrid variants (see, e.g., the review article~\cite{XZ17}), can typically not significantly speed up convergence. This is because they do not sufficiently consider the geometry of the underlying problem. Also, preconditioners for HDG methods such as~\cite{CockburnDGT2013,FabienKMR19,Lu2021,LuRK22a,Lu2023,WidleyMB21} are not suitable for the present application since they require a coarsening strategy that, in the spatial network setting, would change the geometry of the underlying graph. Therefore, the injection operators defined therein cannot be used directly. To overcome these difficulties, this paper employs a preconditioner which is based on the observation that the network can be treated essentially as a continuous object at sufficiently coarse scales. This allows one to place an artificial coarse mesh over the network and use finite element techniques with respect to this mesh to introduce, e.g., a two-level overlapping Schwarz preconditioner similar to~\cite{GoHeMa22}. We prove that the global system matrix is spectrally equivalent to an edge length weighted graph Laplacian (in each component). Under certain network connectivity assumptions, this then allows us to prove the uniform convergence of the corresponding preconditioned conjugate gradient method.

Note that, alternatively, multiscale methods such as the (Super-)Localized Orthogonal Decomposition (cf.~\cite{MaP14,HeP13,HaPe21b,pumslod}) could be used to tackle the problem of large and ill-conditioned linear systems of equations, see~\cite{EGHKM24,HaM22,HMM23}. These methods feature a compression property that allow one to significantly reduce the size of the linear systems of equations to be solved, making them  tractable again, e.g., by classical sparse direct solvers. Also the use of other multiscale methods such as the Multiscale Spectral Generalized Finite Element Method \cite{Babuska2011,Ma2022} or Gambles~\cite{Owh17,OwhS19} seems possible for spatial network models.  For a comprehensive overview of multiscale methods, see also~\cite{MalP20,Peterseim2021}.

The rest of this paper is organized as follows: In \cref{SEC:modelproblem}, we introduce a Timoshenko beam network model, which is then discretized using a HDG method in \cref{SEC:hdg}. An a priori error analysis of the method is given in \cref{SEC:error}, while an efficient preconditioner is presented in \cref{SEC:precond}. The paper concludes with some numerical experiments in \cref{SEC:numexp}. 

\section{Timoshenko beam network model}\label{SEC:modelproblem}
%
\begin{figure}
 \includegraphics[width=.5\textwidth]{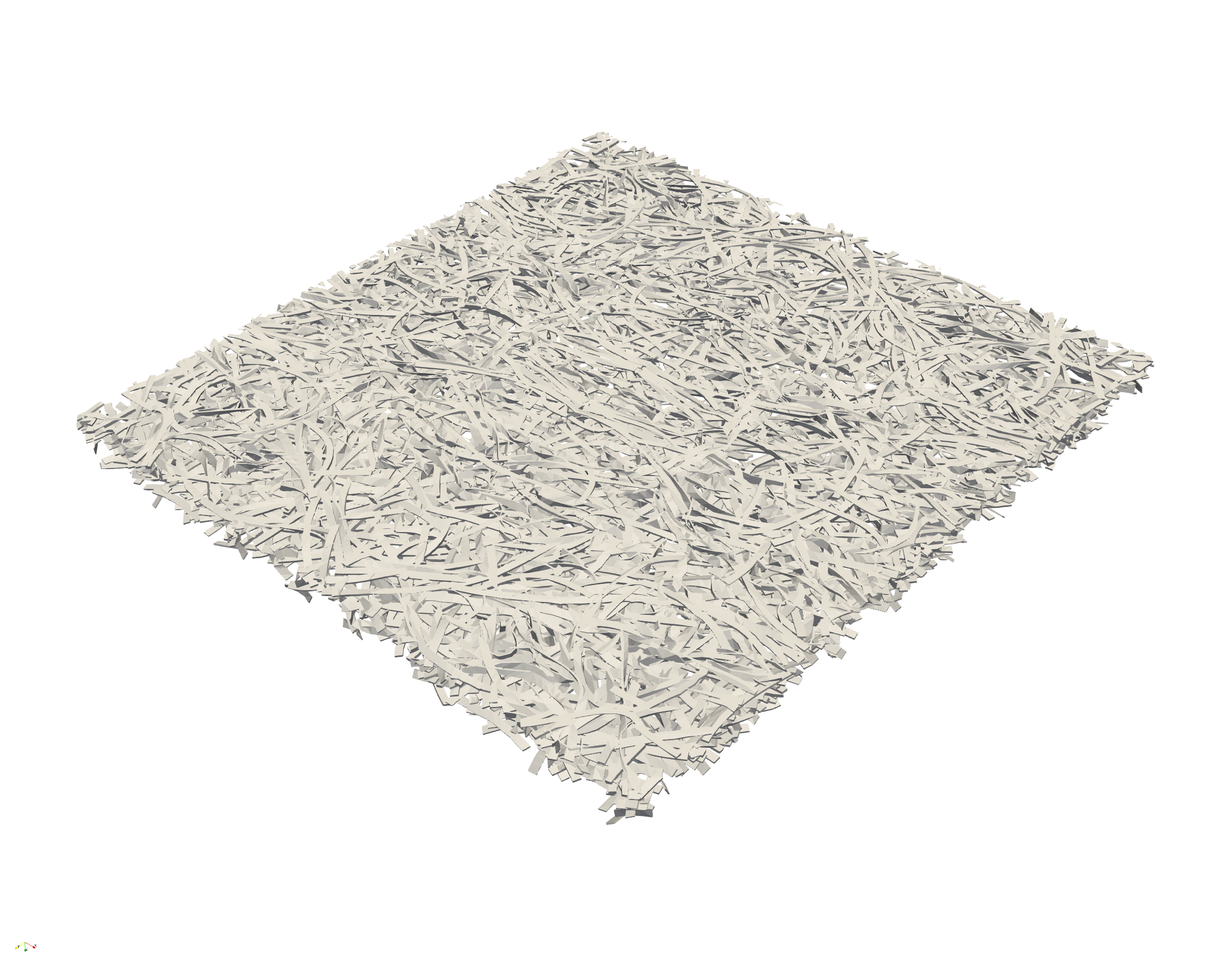}
 \caption{Fiber network model of paper at the millimeter scale.}\label{fig:paper}
\end{figure}
In this section, we derive a Timoshenko beam network model that describes the elastic deformation of fiber-based materials such as paper or cardboard. The beam network model considers a finite number of initially straight beams whose cross sections are constant along the length of the beams. Note that beams are slender objects, i.e., their axial dimension is much larger than their cross sectional diameter. This motivates to represent the network of beams by a graph $\mathcal G = (\mathcal N,\mathcal E)$, where $\setNode = \{ \node_1, \node_2, \dots, \node_K \}$ is a set of zero-dimensional nodes and $\setEdge = \{ \edge_1, \edge_2, \dots. \edge_L\}$ is a set of locally one-dimensional edges. 
We assume that the graph~$\mathcal G$ is connected, which is natural, since otherwise its connected components could be considered separately. 
In the following, a three-dimensional setting is considered, i.e., the nodes and edges of~$\mathcal G$ are contained in~$\mathbb R^3$. Each edge connects two different nodes, and we write $\node \sim \edge$ if the node $\node \in \setNode$ is an endpoint of the edge $\edge \in \setEdge$. In the present paper application, the nodes~$\mathfrak n \in \mathcal N$ model the joints of the beams, whereas the edges $\edge\in \setEdge$ represent the fiber segments connecting two nodes; we refer to \cref{fig:paper} for an illustration of a fiber network model of paper. 
To describe the deformation of the beam network in response to applied forces and moments, we use the Timoshenko beam theory originally developed in~\cite{Timoshenko1921}; see also~\cite{Carrera2011}. In contrast to Euler--Bernoulli beam theory (see, e.g.,~\cite{Bauchau2009}), it can accurately capture shear effects. 

\subsection{Governing equations}\label{SEC:equations}
The derivation of Timoshenko beam theory is based on the laws of linear elasticity. It relies on the assumption that the cross sections of a beam are infinitely rigid in their plane and remain plain after deformation. Furthermore, it is always assumed that the cross sections can rotate independently from the deformation of the centroid line. Under these assumptions, Timoshenko beam theory states for each beam one-dimensional differential equations describing its deformation. The deformation of the beam corresponding to edge $\edge \in \setEdge$ is described by the beam's centerline displacement $\vec u_\edge \colon \edge \to \IR^3$ and its cross section rotation $\vec \re \colon \edge \to \IR^3$. These variables will serve as primal unknowns in the Timoshenko beam equations. The corresponding dual unknowns are $\vec n_\edge \colon \edge \to \IR^3$, the stress resultants from normal and shear forces, and $\vec m_\edge \colon \edge \to \IR^3$, the resultant moment from torsion and bending moments. Denoting by $\node_k$ and~$\node_\ell$ with $k < \ell$ the endpoints of the edge $\edge$, we can define the unit vector aligned with~$\edge$ as ${\vec i_\edge \coloneqq {(\node_\ell - \node_k)}\big/\edgelength}$, where $\edgelength \coloneqq {| \node_\ell - \node_k |}$ denotes the length of $\edge$ and $|\cdot|$ is the Euclidean norm. We further define the (scalar) unit normals of the edge $\edge$ by  $\nu_\edge(\node_k) \coloneqq -1$ and $\nu_\edge(\node_\ell) \coloneqq +1$. 

The properties of the beams are determined by their material and shape parameters, which are encoded into the coefficients~$C_{\vec n}$ and~$C_{\vec m}$; see also \cref{rem:locform}. In the following, we assume that~$C_{\vec n}$ and~$C_{\vec m}$ are symmetric $\mathbb R^{3\times 3}$-valued functions which satisfy uniform lower and upper spectral bounds, i.e., there exist  $0<\alpha_{\vec n},\alpha_{\vec m},\beta_{\vec n},\beta_{\vec m}<\infty$ such~that
\begin{equation}\label{eq:spectralbounds}
 \alpha_{\vec n} |\xi|^2 \leq (C_{\vec n}(x) \xi) \cdot \xi \leq \beta_{\vec n} |\xi|^2,\qquad \alpha_{\vec m} |\xi|^2 \leq (C_{\vec m}(x) \xi) \cdot \xi \leq \beta_{\vec m} |\xi|^2
\end{equation}
holds for all $x \in \edge$ and $\edge \in \setEdge$.

Given the  distributed force $\vec f_\edge \colon \edge \to \IR^3$ and the distributed moment $\vec g_\edge \colon \edge \to \IR^3$, the one-dimensional Timoshenko beam equations corresponding to edge $\edge \in \setEdge$ read
\begin{subequations}\label{EQ:timo_network}
\begin{align}
 - C_{\vec n} ( \partial_x  \vec u_\edge + \vec i_\edge \times \vec \re ) &= \vec n_\edge, & - C_{\vec m} \partial_x \vec \re &= \vec m_\edge, \label{EQ:timo_primal_dual}\\
 \partial_x \vec n_\edge &= \vec f_\edge, & \partial_x \vec m_\edge  + \vec i_\edge \times \vec n_\edge & = \vec g_\edge, \label{EQ:timo_force}
\end{align} 
where $\partial_x$ denotes differentiation with respect to the spatial variable $x$ that varies along~$\edge$ with unit speed and $\times$ denotes the cross product.
We impose continuity and balance conditions at the nodes to couple the one-dimensional differential equations on each edge. The \emph{continuity conditions} enforce that the displacements and rotations at each node coincide. Following the concept of hybridization, this can be reformulated by requiring that for each node $\node \in \setNode$ there exist vectors $\vec \un$, $\vec \rn \in \mathbb R^3$ such that
\begin{equation}\label{EQ:timo_cont_hybrid}
 \vec u_\edge (\node) = \vec \un, \qquad \vec \re (\node) = \vec \rn
\end{equation}
holds for all edges $\edge \in \setEdge$ adjacent to $\node$. The nodewise unknowns $\vec \un$ and $\vec \rn$ will act as Lagrange multipliers in a hybrid formulation of the problem.
We also prescribe Dirichlet boundary conditions at a given set of Dirichlet nodes $\emptyset \neq \setNodeDir \subset \setNode$. More precisely, for each node $\node \in \setNodeDir$ and given Dirichlet data $\und, \rnd \in \mathbb R^3$, we require that 
\begin{equation}\label{EQ:timo_dir}
 \vec \un = \und, \qquad \vec \rn = \rnd.
\end{equation}
The \emph{balance conditions} ensure an equilibrium of forces and moments at non-Dirichlet nodes. This means that for all $\node \in \setNode\setminus \setNodeDir$ and given concentrated forces and moments $\vec f_\node, \vec g_\node \in \mathbb R^3$, which are typically set to zero, it should hold that 
\begin{equation}\label{EQ:timo_balance}
 - \jump{ \vec n_\edge \nu_\edge }_\node = - \vec f_\node, \qquad - \jump{\vec m_\edge \nu_\edge }_\node = - \vec g_\node,
\end{equation}
\end{subequations}
where the summation operator $\jump{\cdot}_\node$ sums over all values that are attained at $\node$. This notation is a generalization of the jump operator in the context of discontinuous Galerkin methods; see, e.g.,~\cite[Eq.~(2.10)]{RuppGK22}. The multiplication of equation~\cref{EQ:timo_balance} with $-1$ is purely artificial at this point, but it will allow us later to deal with positive definite instead of negative definite operators. In the engineering literature, nodes subject to conditions \cref{EQ:timo_cont_hybrid,EQ:timo_balance} are typically called \emph{rigid joints}, and elastic structures with rigid joints are sometimes called \emph{frames}. Problem~\cref{EQ:timo_network} defines a frame of Timoshenko beams. 

Remembering that $\mathcal G$ is connected, there is a unique solution to \cref{EQ:timo_network}, since there is at least one Dirichlet node.

\begin{remark}[Local formulation]\label{rem:locform}
 In the literature, problem~\cref{EQ:timo_network} is typically formulated in local coordinates, which requires a change of basis between local and global coordinates. As local basis corresponding to edge  $\edge \in \setEdge$ we consider $\{\vec i_\edge,\vec j_\edge,\vec k_\edge\}$, where $\vec i_\edge$ is as above, and~$\vec j_\edge$ and~$\vec k_\edge$ are chosen as the principal axes of inertia of the beam's cross section. As global basis the canonical basis of~$\mathbb R^3$ is chosen, which we denote by $\{\hat{\vec i}, \hat{\vec j}, \hat{\vec k}\}$. We assume that  $\{\vec i_\edge,\vec j_\edge,\vec k_\edge\}$ is right-handed, i.e., it holds that $\vec i_\edge = \vec j_\edge \times \vec k_\edge $.  The change-of-basis matrix for the beam $\edge$ is then given by $\vec T_\edge \coloneqq (\vec i_\edge,\vec j_\edge,\vec k_\edge) \in \mathbb R^{3\times 3}$. Since the basis is right-handed and orthonormal, it holds that $\det(\vec T_\edge) = 1$. Multiplying \eqref{EQ:timo_primal_dual} and \eqref{EQ:timo_force} from the left by~$\vec T_\edge^\top$ and using that $\vec T_\edge$ is orthogonal, we arrive at the local problem
 \begin{align}\label{EQ:timolocal}
  - \hat C_{\vec n} ( \partial_{\hat x}  \hat{\vec u}_\edge + \hat{\vec i} \times \hat{\vec r}_\edge ) &= \hat{\vec n}_\edge, & - \hat C_{\vec m} \partial_{\hat x} \hat{\vec r}_\edge &= \hat{\vec m}_\edge,\\
  \partial_{\hat x} \hat{\vec n}_\edge &= \hat{\vec f}_\edge, &  \partial_{\hat x} \hat{\vec m}_\edge+ \hat{\vec i} \times \hat{\vec n}_\edge & = \hat{\vec g}_\edge,
 \end{align} 
 where local variables are indicated by hats. This problem is posed on $[0, \edgelength] \times \{ 0 \}^2$. Note that due to the choice of local basis, the coefficient matrices $\hat C_{\vec n}(x) = \vec T_\edge^\top C_{\vec n}(x) \vec T_\edge$ and $\hat C_{\vec m}(x)= \vec T_\edge^\top C_{\vec m}(x)\vec T_\edge$ are diagonal for all $x \in \edge$  and the components of $\hat{\vec u}_\edge $ and $\hat{\vec r}_\edge$ decouple.
 More precisely, we have the characterization
 \begin{equation*}
  \hat C_{\vec n} = \operatorname{diag}(EA, kGA, kGA),\qquad \hat C_{\vec m} = \operatorname{diag}(GI_t, EI_2, EI_3),
 \end{equation*}
 where $E$ is the elastic modulus, $G$ is the shear modulus, $A$ is the cross section area, $kA$ is the corrected shear area, $I_2$ and $I_3$ are the second moments of inertia of the cross section and $I_t$ is the polar moment of the cross section. Note that the local formulation~\cref{EQ:timolocal} is essential for implementing Timoshenko beam networks. 
\end{remark}

\section{HDG discretization}\label{SEC:hdg}
%
This section derives an HDG method for the numerical solution of the Timoshenko beam network model introduced in the previous section. To simplify the presentation, the following remark introduces a notation to hide constants in estimates independent of the edge length $\edgelength$ and the stabilization parameter $\tau_\edge$ (to be introduced later).

\begin{remark}[Tilde notation]
	If it  holds that $a \leq C b$, where $C>0$ is a constant that may depend on the coefficients $C_{\vec n}$ and $C_{\vec m}$, the seminorms of the solutions $\ue, \re, \vec n_\edge,$ and~$\vec m_\edge$, but is independent of the edge length $\edgelength$ and the stabilization parameter $\tau_\edge$, we may write $a \lesssim b$ to hide the constant. Similarly, we may write $b\gtrsim a$ for $a \geq C b$.
\end{remark}

\subsection{Hybrid dual mixed formulation}\label{SEC:hybrid_formulation}
In the following, we present the hybrid dual-mixed formulation of the Timoshenko beam network model, which will serve as the starting point for deriving the HDG method. It is based on the function spaces ${V^\edge_{\vec u} \coloneqq  (L^2(\edge))^3}$  and $V^\edge_{\vec n} \coloneqq (H^1(\edge))^3$ defined locally on each edge $\edge \in \setEdge$. Let us also introduce the $L^2(\edge)$-inner product and a lower-dimensional version of it acting on the nodes at the endpoints of $\edge$ for all functions $\vec v, \vec w \in V_{\vec u}^\edge$ and $\vec p, \vec q \in V^\edge_{\vec n}$ by
\begin{align*}
 (\vec v,\vec w)_\edge \coloneqq \int_\edge \vec v \cdot \vec w \ds,\qquad \langle \vec p,\vec q\rangle_\edge \coloneqq \sum_{\node \sim \edge} \vec p(\node) \cdot \vec q(\node),
\end{align*}
where the point evaluation is meaningful due to the continuous embedding ${H^1(\edge)\hookrightarrow \mathcal C^0(\overline{\edge})}$ which holds in one spatial dimension. 
Note that in an abuse of notation, we will also allow functions defined only on the nodes of the networks to be plugged into $\langle \cdot,\cdot\rangle_\edge$, such as functions involving the unit normal~$\nu_\edge$, but also $\un$ and $\rn$. 

Given the  data
	\begin{itemize}
		\item $\vec f_\edge, \vec g_\edge \in V_{\vec u}^\edge$ for all edges $\edge \in \setEdge$ ,
		\item $\vec f_\node, \vec g_\node \in \IR^3$ for all free nodes $\node \in \setNode \setminus \setNodeDir$,
		\item $\und, \rnd\in \mathbb R^3$ for all Dirichlet nodes $\node \in \setNodeDir$,
	\end{itemize}
the hybrid dual mixed formulation of problem \cref{EQ:timo_network} seeks
\begin{itemize}
 \item $\vec u_\edge, \vec \re \in V^\edge_{\vec u}$ for all edges $\edge \in \setEdge$, 
 \item $\vec n_\edge, \vec m_\edge \in V^\edge_{\vec n}$ for all edges $\edge \in \setEdge$, 
 \item $\vec \un, \vec \rn \in \IR^3$ for all nodes $\node \in \setNode$,
\end{itemize}
we refer to as primal, dual, and hybrid unknowns, respectively. The hybrid unknowns should satisfy for all $\node \in \setNodeDir$ the Dirichlet boundary conditions~\cref{EQ:timo_dir}, while the dual unknowns should satisfy for all $\node \in \setNode \setminus \setNodeDir$ the balance conditions \cref{EQ:timo_balance}. Moreover, for all edges $\edge \in \setEdge$ the following equations should hold for all $\vec p, \vec q \in V^\edge_{\vec n}$ and  $\vec v, \vec w \in V^\edge_{\vec u}$:
\begin{equation}\label{EQ:local_solver}
 \arraycolsep=1.4pt\def\arraystretch{2.2}
 \begin{array}{clclclclcl}
  -&(C^{-1}_{\vec n} \vec n_\edge, \vec p)_\edge &&& + & (\vec u_\edge, \partial_x \vec p)_\edge & - & (\vec i_\edge \times \vec \re, \vec p)_\edge & = &
  \langle\un,\vec p \nu_ \edge\rangle_\edge
   , \\[-.75em]
  && - & (C^{-1}_{\vec m} \vec m_\edge, \vec q)_\edge &&& + & (\vec \re, \partial_x \vec q)_\edge & = & 
  \langle \rn,\vec q\nu_\edge\rangle_\edge
  , \\[-.75em]
  &(\partial_x \vec n_\edge, \vec v)_\edge &&&&&&& = & (\vec f_\edge, \vec v)_\edge, \\[-.75em]
  &(\vec i_\edge \times \vec n_\edge, \vec w)_\edge & + & (\partial_x \vec m_\edge, \vec w)_\edge &&&&& = & (\vec g_\edge, \vec w)_\edge.
 \end{array}
\end{equation}
Note that from the first and second equations in~\cref{EQ:local_solver}, one can immediately conclude that $\vec u_\edge, \vec \re \in (H^1(\edge))^3$, even though we only seeked functions in $(L^2(\edge))^3$.

In the following, we will regularly view problem~\cref{EQ:local_solver} as a local solver that maps the Dirichlet data and distributed source terms corresponding to an edge to the solution on that edge. The following lemma proves the well-posedness of this local solver.

\begin{lemma}[Well-posedness of local solver]\label{LEM:loc_ex}
 Consider the edge $\edge \in \setEdge$ and let the boundary data $\vec \un, \vec \rn \in \IR^3$ for the two nodes $\node$ at the endpoints of $\edge$ as well as the source terms $\vec f_\edge, \vec g_\edge \in V_{\vec u}^\edge$ be given. Then, the local solver defined by problem \cref{EQ:local_solver} has a unique solution comprised of $\vec u_\edge, \vec \re \in V^\edge_{\vec u}$ and $\vec n_\edge, \vec m_\edge \in V^\edge_{\vec n}$. 
\end{lemma}
\begin{proof}
 As a first step, we reformulate the local solver \cref{EQ:local_solver} as a saddle-point problem, which fits into the framework of classical inf-sup theory, cf.~\cite{Boffi2013}. For this, we multiply the whole system by $-1$ and  define the following bilinear and linear forms:
 \begin{align*}
  a((\vec n_\edge, \vec m_\edge), (\vec p, \vec q)) &\coloneqq  (C^{-1}_{\vec n} \vec n_\edge, \vec p)_\edge + (C^{-1}_{\vec m} \vec m_\edge, \vec q)_\edge,\\
  b((\vec n_\edge, \vec m_\edge),(\vec v, \vec w))&\coloneqq -(\partial_x \vec n_\edge, \vec v)_\edge - (\vec i_\edge \times \vec n_\edge, \vec w)_\edge - (\partial_x \vec m_\edge, \vec w)_\edge.
 \end{align*}
 Using the properties of the cross product, we can then write \cref{EQ:local_solver} as the saddle-point problem, which seeks $(\vec u_\edge,\vec r_\edge) \in V_{\vec u}^\edge\times V_{\vec u}^\edge$ and $(\vec n_\edge,\vec m_\edge) \in V_{\vec n}^\edge\times V_{\vec n}^\edge$ such that \vspace{-1.5ex}
 	\begin{equation}\label{EQ:bil_forms}
 		\arraycolsep=1.4pt\def\arraystretch{2.2}
 		\begin{array}{lclcl}
 			a((\vec n_\edge, \vec m_\edge), (\vec p, \vec q)) & + & b((\vec p, \vec q),(\vec u_\edge, \vec \re))) & = & -\langle\un,\vec p\nu_\edge\rangle_\edge-\langle\rn,\vec q \nu_\edge\rangle_\edge, \\[-.75em]
 			b((\vec n_\edge, \vec m_\edge),(\vec v, \vec w)) &&& = & - (\vec f_\edge, \vec v)_\edge - (\vec g_\edge, \vec w)_\edge
 		\end{array}
 	\end{equation}
 	holds for all $(\vec v,\vec w) \in V_{\vec u}^\edge\times V_{\vec u}^\edge$ and $(\vec p,\vec w) \in V_{\vec n}^\edge\times V_{\vec n}^\edge$.

 To prove the well-posedness of \cref{EQ:bil_forms}, we verify the two BNB conditions, i.e., the coercivity of $a$ on the kernel of $b$ and the inf-sup stability of $b$. These conditions involve the tuple spaces $V_{\vec u}^\edge\times V_{\vec u}^\edge$ and $V_{\vec n}^\edge\times V_{\vec n}^\edge$, whose canonical norms we denote by 
 \begin{align*}
  \|(\vec n_\edge,\vec m_\edge)\|_{V_{\vec n}^\edge\times V_{\vec n}^\edge}^2 &\coloneqq \|\vec n_\edge\|_\edge^2 + \|\partial_x \vec n_\edge\|_\edge^2 + \|\vec m_\edge\|_\edge^2 + \|\partial_x \vec m_\edge\|_\edge^2,\\
  \|(\vec v,\vec w)\|_{V_{\vec u}^\edge\times V_{\vec u}^\edge}^2 &\coloneqq  \|\vec v\|_\edge^2 + \|\vec w\|_\edge^2.
 \end{align*}
 For proving the first BNB condition, which is the coercivity of $a$ on $\operatorname{ker}(b)$, we note that $(\vec n_\edge,\vec m_\edge) \in \operatorname{ker}(b)$ implies that $\partial_x \vec n_\edge = 0$ and $\partial_x \vec m = -\vec i_\edge \times \vec n_\edge$. This, together with the uniform spectral bounds for~$C_{\vec n}$ and $C_{\vec m}$, cf.~\cref{eq:spectralbounds}, yields for any $(\vec n_\edge,\vec m_\edge) \in \operatorname{ker}(b)$ that
 \begin{align*}
  a((\vec n_\edge, \vec m_\edge),(\vec n_\edge, \vec m_\edge)) &\gtrsim \|\vec n_\edge\|_\edge^2 + \|\partial_x \vec n_\edge\|_\edge^2 + \tfrac12\|\vec m_\edge\|_\edge^2  + \tfrac12\|\partial_x \vec m_\edge\|_\edge^2 - \tfrac12\|\vec i_\edge \times \vec n_\edge\|_\edge^2\\
  &\geq \tfrac12\|(\vec n_\edge,\vec m_\edge)\|_{V_{\vec n}^\edge\times V_{\vec n}^\edge}^2,
 \end{align*}
 which is the desired coercivity property. 

To prove the second BNB condition, which is the inf-sup stability of $b$, we will choose for any given tuple $(\vec v,\vec w) \in V_{\vec u}^\edge\times V_{\vec u}^\edge$ a tuple~$(\vec n_\edge, \vec m_\edge) \in V_{\vec n}^\edge \times V_{\vec n}^\edge$ whose components $\vec n_\edge$ and $\vec m_\edge$ are defined as the antiderivatives of $\vec v$ and~$\vec w$ along the edge $\edge$,  respectively. The integration constant of the antiderivatives is set to zero. This yields for any $(\vec v,\vec w) \in V_{\vec u}^\edge\times V_{\vec u}^\edge$ that
 \begin{align*}
  \sup_{(\vec n_\edge,\vec m_\edge)\neq 0}\frac{b((\vec n_\edge,\vec m_\edge),(\vec v,\vec w))}{\|(\vec n_\edge,\vec m_\edge)\|_{V_{\vec n}^\edge\times V_{\vec n}^\edge}} \geq \big(1+\tfrac{4\edgelength^2}{\pi^2}\big)^{-1/2}\;\frac{(\vec v,\vec v)_\edge + (\vec w,\vec w)_\edge  + (\vec i_\edge\times \vec n_\edge,\vec w)_\edge}{\|(\vec v,\vec w)\|_{V_{\vec u}^\edge\times V_{\vec u}^\edge}},
 \end{align*}
 where we used that $\|\vec n_\edge\|_\edge \leq \tfrac{2\edgelength}{\pi}\|\vec v\|_\edge$ and $\|\vec m_\edge\|_\edge \leq \tfrac{2\edgelength}{\pi}\|\vec w\|_\edge$. The latter constant can be derived by explicitly computing the first eigenvalue of the Laplacian eigenvalue problem on the edge $\edge$ subject to mixed boundary conditions. Note that the estimates also imply that $(\vec i_\edge\times \vec n_\edge,\vec w)_\edge \leq \|\vec n_\edge\|_\edge \|\vec w\|_\edge \leq \tfrac{2\edgelength}{\pi} \|\vec v\|_\edge\|\vec w\|_\edge$, which can be used to estimate the numerator on the right-hand side of the previous inequality as 
 \begin{align*}
  (\vec v,\vec v)_\edge + (\vec w,\vec w)_\edge  + (\vec i_\edge\times \vec n_\edge,\vec w)_\edge \geq \big(1-\tfrac{\edgelength}\pi\big)\|(\vec v,\vec w)\|_{V_{\vec u}^\edge\times V_{\vec u}^\edge}^2.
 \end{align*}
 Thus, the bilinear form $b$ is inf-sup stable if the edge length satisfies $\edgelength <\pi$. Note that this limitation can be easily overcome by rescaling the network, so it is not assumed in the following. The assertion can be concluded using the classical inf-sup theory.
\end{proof}

In its above form, the hybrid dual mixed formulation involves primal, dual, and hybrid unknowns.  To derive an equivalent condensed problem involving only hybrid unknowns, we use an alternative interpretation of the local solver defined in~\cref{EQ:local_solver}. More precisely, for each edge $\edge \in \setEdge$, the local solver defines the following mappings:
\begin{equation}\label{EQ:locmapping}
 \begin{aligned}
 	  \vec U_\edge\colon (\vec \un, \vec \rn, \vec f_\edge, \vec g_\edge) &\mapsto \vec u_\edge, \qquad &\vec R_\edge\colon (\vec \un, \vec \rn, \vec f_\edge, \vec g_\edge) &\mapsto \vec \re,\\
  \vec N_\edge\colon (\vec \un, \vec \rn, \vec f_\edge, \vec g_\edge) &\mapsto \vec n_\edge, \qquad &\vec M_\edge\colon (\vec \un, \vec \rn, \vec f_\edge, \vec g_\edge) &\mapsto \vec m_\edge.
 \end{aligned}
\end{equation}
Note that we abuse the notation and require that each of the first two arguments contains the two values corresponding to the nodes at the endpoints of $\edge$. In the case of homogeneous source terms, i.e., $\vec f_\edge = \vec g_\edge = 0$, we will omit the last two arguments.

Before we state the condensed problem, we need to introduce some notation. We denote the trial functions which are set to zero at Dirichlet nodes by $\vec \lambda, \vec \phi \colon \setNode \to \mathbb R^3$. Writing $\vec \lambda_\node$ and $\vec \phi_\node$ for the values of $\vec \lambda$ and $\vec \phi$ at the node $\node$, respectively, this means that $\vec \lambda_\node = \vec \phi_\node = 0$ for all $\node \in\setNodeDir$.  Similarly, the test functions denoted by $\vec \mu,\vec \psi\colon \setNode\to \mathbb R^3$ should also satisfy that $\vec \mu_\node= \vec \psi_\node=0$ for all $\node \in \setNodeDir$. Denoting the function space the trial and test functions are contained in by $V_{\vec \lambda}$, the condensed problem seeks $(\vec \lambda,\vec \phi) \in V_{\vec \lambda}\times V_{\vec \lambda}$ such that  
\begin{equation}\label{EQ:hybrid_problem}
 A((\vec \lambda, \vec \phi),(\vec \mu, \vec \psi)) = F((\vec \mu, \vec \psi))
\end{equation}
holds for all $(\vec \mu,\vec \psi) \in V_{\vec \lambda}\times V_{\vec \lambda}$. 
The bilinear form $A$ of the condensed problem encodes the balance of forces and moments conditions for the operator-harmonic (homogeneous source terms) extension of the nodal data $(\vec \lambda,\vec \phi)$, cf.~\cref{EQ:timo_balance}, and is given by
\begin{align}\label{EQ:bilinearforms}
 A((\vec \lambda, \vec \phi),(\vec \mu, \vec \psi)) \coloneqq  -\sum_{\node \in \setNode \setminus \setNodeDir}  & \left[\jump{ \vec N_\edge(\vec \lambda_\node, \vec \phi_\node)\nu_\edge}_\node\cdot  \vec \mu_\node + \jump{ \vec M_\edge(\vec \lambda_\node, \vec \phi_\node)\nu_\edge }_\node\cdot \vec \psi_\node\right].
\end{align}
Inhomogeneous Dirichlet boundary data and source terms, as well as prescribed values for the balance conditions, are incorporated by the bilinear form $F$ defined by
\begin{align*}
 \begin{split}
  F((\vec \mu, \vec \psi)) \coloneqq\sum_{\node \in \setNode \setminus \setNodeDir} &\big[\hspace{-.4ex}\left(\jump{ \vec N_\edge(\und, \rnd, \vec f_\edge, \vec g_\edge)\nu_\edge}_\node - \vec f_\node\right)\cdot  \vec \mu_\node \\[-2ex]
  &\quad+ \left(\jump{ \vec M_\edge(\und, \rnd, \vec f_\edge, \vec g_\edge)\nu_\edge}_\node - \vec g_\node\right) \cdot \vec \psi_\node\big],
 \end{split}
\end{align*}
where we set $\und = \rnd = 0$ for all $\node \in \setNode\setminus \setNodeDir$. After having computed $(\vec \lambda,\vec \phi)$ as the solution to \cref{EQ:hybrid_problem}, the solution to the hybrid dual mixed formulation can be retrieved as follows: The displacements and rotations at the nodes can be obtained as $\un = \vec \lambda_\node + \und$ and $\rn = \vec \phi_\node + \rnd$ for all $\node \in \setNode$ and the unknowns $\ue, \re,\vec n_\edge,\vec m_\edge$ can for all $\edge \in \setEdge$ be obtained from $\un$ and $\rn$ using the local mappings defined in \cref{EQ:locmapping}, i.e.,
\begin{align}\label{EQ:retrievsolution}
 \begin{aligned}
 	  \ue &= \vec U_\edge(\un, \rn, \vec f_\edge, \vec g_\edge),\qquad &\re &= \vec R_\edge(\un, \rn, \vec f_\edge, \vec g_\edge),
\\
  \vec n_\edge &= \vec N_\edge(\un, \rn, \vec f_\edge, \vec g_\edge),\qquad &\vec m_\edge &= \vec M_\edge(\un, \rn, \vec f_\edge, \vec g_\edge).
 \end{aligned}
\end{align} 

The following lemma shows that the bilinear form $A$ is symmetric positive definite, proving the well-posedness of the condensed problem.

\begin{lemma}[Properties of condensed problem]\label{LEM:characterize_bil}
	The bilinear form $A$ can be equivalently rewritten as the following symmetric expression
 \begin{align}\label{eq:identitysym}
  \begin{split}
   A((\vec \lambda, \vec \phi),(\vec \mu, \vec \psi)) = \sum_{\edge \in \setEdge}& \big[(C^{-1}_{\vec n} \vec N_\edge(\vec \lambda_\node, \vec \phi_\node), \vec N_\edge(\vec \mu_\node, \vec \psi_\node))_\edge\\[-2ex]
   &\quad + (C^{-1}_{\vec m} \vec M_\edge(\vec \lambda_\node, \vec \phi_\node), \vec M_\edge(\vec \mu_\node, \vec \psi_\node))_\edge\big].
  \end{split}
 \end{align}
This implies that the bilinear form $A$ is symmetric and positive definite, and hence the well-posedness of problem \cref{EQ:hybrid_problem}.
\end{lemma}
\begin{proof}
 In this proof, we will abbreviate $\vec N_\edge(\vec \lambda_\node,\vec \phi_\node)$ and $\vec N_\edge(\vec \mu_\node,\vec \psi_\node)$ by $\vec N_\edge^1$ and $\vec N_\edge^2$, respectively, and analogous abbreviations will be used for the other local mappings. Testing the first and second equations of system~\cref{EQ:local_solver} with $\Ndisc^1$ and $\Ndisc^2$, respectively and using the symmetry of $C_{\vec n}$ and $C_{\vec m}$, we obtain that
 \begin{align}
  - \langle \vec \mu_\node,\vec N_\edge^1\nu_\edge\rangle_\edge &= (C_{\vec n}^{-1} \vec N_\edge^1,\vec N_\edge^2)_\edge -(\vec U_\edge^2,\partial_x \vec N_\edge ^1)_\edge
  + (\vec i_\edge \times \vec R_\edge^2,\vec N_\edge^1)_\edge,\label{eq:sym1}\\
  - \langle\vec \psi_\node,\vec M_\edge^1\nu_\edge\rangle_\edge &=(C_{\vec m}^{-1} \vec M_\edge^1,\vec M_\edge^2)_\edge -(\vec R_\edge^2,\partial_x \vec M_\edge ^1)_\edge. \label{eq:sym2}
 \end{align}
 To rewrite the latter expressions, we test the third and fourth equations of system~\cref{EQ:local_solver} with $\Udisc^2$ and $\Rdisc^2$, respectively, and note that $\vec f_\edge = \vec g_\edge = 0$. Inserting the resulting equations into \cref{eq:sym1,eq:sym2} and summing up yields the identity
 \begin{align*}
  -\langle \vec \mu_\node,\vec N_\edge^1\nu_\edge\rangle_\edge - \langle\vec \psi_\node,\vec M_\edge^1\nu_\edge\rangle_\edge= (C_{\vec n}^{-1} \vec N_\edge^1,\vec N_\edge^2)_\edge + (C_{\vec m}^{-1} \vec M_\edge^1,\vec M_\edge^2)_\edge.
 \end{align*}
 Equation~\cref{eq:identitysym} then follows after summing over all edges and using the definition of $\jump{\cdot}(\node)$.

 To prove the positive definiteness of the bilinear form $A$, it only remains to show that $A((\vec \lambda,\vec \phi),(\vec \lambda,\vec \phi)) = 0$ implies that $(\vec \lambda,\vec \phi) = 0$. However, this is a consequence of the well-posedness of the local solver \cref{EQ:local_solver} proved in \cref{LEM:loc_ex}. The positive definiteness implies the invertibility of the operator $A$ and thus the well-posedness of problem~\cref{EQ:hybrid_problem}.
\end{proof}

\subsection{Discretization}\label{SEC:discretization}

The HDG discretization for the numerical solution of the Timoshenko beam network model \cref{EQ:timo_network} is derived by discretizing the hybrid dual mixed formulation. Note that only the local solvers defined in~\cref{EQ:local_solver} must be discretized since all other aspects are already finite-dimensional. For the discretization of the local solver corresponding to the edge $\edge \in\setEdge$ we use the space of polynomials defined on $\edge$ of degree at most~$p$, where $p \in \mathbb N$ is a fixed number. This space and the corresponding space of componentwise polynomials are denoted by $\mathds P_p(\edge)$ and $V_p^\edge \coloneqq (\mathds P_p(\edge))^3$, respectively.

Given the data as in \cref{SEC:hybrid_formulation}, the HDG method seeks
\begin{itemize}
 \item $\uedisc, \redisc \in V_p^\edge$ for all edges $\edge \in \setEdge$, 
 \item $\ndisc, \mdisc \in V_p^\edge$ for all edges $\edge \in \setEdge$, 
 \item $\undisc, \rndisc \in \IR^3$ for all  nodes $\node \in \setNode $,
\end{itemize}
where the hybrid HDG unknowns should satisfy for all $\node \in \setNodeDir$ the Dirichlet boundary conditions~\cref{EQ:timo_dir}. 
We further require balance conditions for the numerical fluxes, which will be the discrete counterpart to the balance conditions~\cref{EQ:timo_balance}.
Given a stabilization parameter $\tau_\edge >0$, the numerical fluxes defined at the network nodes additionally incorporate the $\tau_\edge$-weighted difference between the primal and hybrid HDG unknowns to increase the method's stability. The stabilization parameter should always satisfy $\tau_\edge \edgelength \lesssim 1$. The balance of numerical fluxes then reads for all $\node \in \setNode\setminus \setNodeDir$ as
\begin{equation}
	\label{EQ:timo_balance_disc}
 -\jump{\bar{\vec n}_\edge + \tau_\edge (\bar{\vec u}_\edge - \bar{\vec u}_\node)}_n =- f_\node,\qquad -\jump{\bar{\vec m}_\edge + \tau_\edge (\bar{\vec r}_\edge - \bar{\vec r}_\node)}_n = -g_\node.
\end{equation}
Moreover, for all edges $\edge \in \setEdge$ we require that the HDG unknowns are connected by the following local equations, which should hold for all $\bar{\vec p}, \bar{\vec q}, \bar{\vec v}, \bar{\vec w} \in V_p^\edge$:
\begin{equation}\label{EQ:discrete_solver}
 \arraycolsep=1.375pt\def\arraystretch{2.2}
 \begin{array}{clclclclcl}
  -&(C^{-1}_{\vec n} \ndisc, \bar{\vec p})_\edge &&& + & (\uedisc, \partial_x \bar{\vec p})_\edge & - & (\vec i_\edge \times \redisc, \vec p)_\edge & = & 
  \langle \undisc ,\bar{\vec p} \nu_\edge\rangle_\edge
  , \\[-.75em]
  && - & (C^{-1}_{\vec m} \mdisc, \bar{\vec q})_\edge &&& + & (\redisc, \partial_x \bar{\vec q})_\edge & = &
   \langle \rndisc,\bar{\vec q} \nu_\edge\rangle_\edge
   ,  \\[-.75em]
  &(\partial_x \ndisc, \bar{\vec v})_\edge &&& + &  \tau_\edge
  \langle\uedisc,\bar{\vec v}\rangle_\edge
   &&& = & (\vec f_\edge, \bar{\vec v})_\edge + \tau_\edge
   \langle\undisc,\bar{\vec v}\rangle_\edge
   , \\[-.75em]
  &(\vec i_\edge \times \ndisc, \bar{\vec w})_\edge & + & (\partial_x \mdisc, \bar{\vec w})_\edge &&& + & \tau_\edge \langle\redisc, \bar{\vec w}\rangle_\edge & = & (\vec g_\edge, \bar{\vec w})_\edge + \tau_\edge \langle\rndisc,\bar{\vec w}\rangle_\edge
  .
 \end{array}
\end{equation}
Compared to \cref{EQ:local_solver}, this system additionally incorporates $\tau_\edge$-weighted stabilization terms, strengthening its diagonal.  This is necessary because our choice of function spaces might render the proof of \cref{LEM:loc_ex} invalid, which relied on an inf-sup condition between the spaces $V_{\vec u}^\edge$ and $V^\edge_{\vec n}$. Problem~\cref{EQ:discrete_solver} defines a discrete version of the local solver from~\cref{EQ:local_solver}, which maps the Dirichlet data and distributed source terms corresponding to an edge to the discrete solution on that edge. The following lemma proves this discrete local solver is well-posed for any positive stabilization parameter. 

\begin{lemma}[Well-posedness of discretized local solver]\label{LEM:disc_solver}
 Consider the edge $\edge \in \setEdge$ and let the boundary data $\undisc, \rndisc \in \IR^3$ for the two nodes $\node$ at the endpoints of $\edge$ as well as the source terms $\vec f_\edge, \vec g_\edge \in V_{\vec u}^\edge$ be given. Then, for any stabilization parameters $\tau_\edge>0$, the discretized local solver defined by~\cref{EQ:discrete_solver} has a unique solution comprised of ${\uedisc, \redisc, \ndisc,\mdisc\in V^\edge_p}$.
\end{lemma}
\begin{proof}
 Since~\cref{EQ:discrete_solver} induces a square linear system of equations, it is sufficient to prove that $\undisc = \rndisc = 0$ and $\vec f_\edge = \vec g_\edge = 0$ implies that $\uedisc = \redisc = \ndisc = \mdisc = 0$. To this end, we test the individual equations of system~\cref{EQ:discrete_solver} with the test functions $-\ndisc$, $-\mdisc$, $\uedisc$, and~$\redisc$, respectively. Summing up the resulting equations then yields that
 \begin{equation}
 	\label{eq:zerobdrcond}
  (C^{-1}_{\vec n} \ndisc, \ndisc)_\edge + (C^{-1}_{\vec m} \mdisc, \mdisc)_\edge + \tau_\edge \langle\uedisc, \uedisc\rangle_\edge + \tau_\edge \langle\redisc, \redisc\rangle_\edge = 0.
 \end{equation}
 By the positive definiteness of the matrices $C_{\vec n}$ and $C_{\vec m}$, cf.~\cref{eq:spectralbounds}, this implies that $\ndisc = \mdisc = 0$. Integrating the second equation in~\cref{EQ:discrete_solver} by parts and exploiting that~$\redisc$ satisfies zero boundary conditions, cf.~\cref{eq:zerobdrcond}, we observe that $\redisc = 0$. Finally, by the very same argument, the first equation of~\cref{EQ:discrete_solver} reveals that $\uedisc = 0$, which concludes the proof.
\end{proof}

To establish a condensed formulation of the HDG method similar to \cref{EQ:hybrid_problem}, we introduce for all edges $\edge \in \setEdge$ discrete analogs to the local mappings from~\cref{EQ:locmapping}. These are defined in terms of the discrete local solver~\cref{EQ:discrete_solver} and are denoted by
\begin{equation}\label{EQ:disclocmapping}
 \begin{aligned}
\Udisc\colon (\undisc, \rndisc, \vec f_\edge, \vec g_\edge) &\mapsto \uedisc, \qquad &\Rdisc\colon (\undisc, \rndisc, \vec f_\edge, \vec g_\edge) &\mapsto \vec \redisc, \\
    \Ndisc\colon (\undisc, \rndisc, \vec f_\edge, \vec g_\edge) &\mapsto \ndisc, \qquad &\Mdisc\colon (\undisc, \rndisc, \vec f_\edge, \vec g_\edge) &\mapsto \mdisc.
 \end{aligned}
\end{equation}
Recalling that $V_{\vec \lambda}$ denotes the space of $\mathbb R^3$-valued functions defined on $\setNode$ with homogeneous Dirichlet boundary conditions at the nodes in the set $\setNodeDir$, the condensed formulation of the HDG method seeks the tuple $(\lambdadisc,\phidisc) \in V_{\vec \lambda}\times V_{\vec \lambda}$ such that  
\begin{equation}\label{EQ:hybrid_problem_disc}
 \bar A((\lambdadisc, \phidisc),(\mudisc, \psidisc)) = \bar F((\mudisc, \psidisc))
\end{equation}
holds for all $(\mudisc,\psidisc) \in V_{\vec \lambda}\times V_{\vec \lambda}$. The bilinear form $\bar A$ encodes the balance conditions of the numerical fluxes at the free nodes of the network and is defined as
\begin{align}\label{EQ:disc_bil_form}
	\begin{split}
		\bar A((\lambdadisc, \phidisc),(\mudisc, \psidisc) 
		\coloneqq- \sum_{\node \in \setNode \setminus \setNodeDir}  &\big[\big(\jump{ \Ndisc(\lambdadisc_\node, \phidisc_\node)\nu_\edge + \tau_\edge ( \Udisc(\lambdadisc_\node, \phidisc_\node) - \lambdadisc_\node) }_\node \cdot \mudisc_\node \\[-2ex]
		& \quad+\jump{  \Mdisc(\lambdadisc_\node, \phidisc_\node)\nu_\edge + \tau_\edge ( \Rdisc(\lambdadisc_\node, \phidisc_\node) - \phidisc_\node) }_\node \cdot \psidisc_\node\big)\big].
			\end{split}
	\end{align}
The linear form $\bar F$, which is used to incorporate inhomogeneous data, is defined as
\begin{align*}
 \begin{split}
  \bar F((\mudisc, \psidisc)) \coloneqq\sum_{\node \in \setNode \setminus \setNodeDir} &\big[\big(\jump{ \bar{\vec N}_\edge(\und, \rnd, \vec f_\edge, \vec g_\edge)\nu_\edge + \tau_\edge \bar{\vec U}_\edge(\und, \rnd, \vec f_\edge, \vec g_\edge) }_\node  - \vec f_\node\big) \cdot \mudisc_\node \\[-2ex]
  &\quad+ \big(\jump{ \Mdisc(\und, \rnd, \vec f_\edge, \vec g_\edge)\nu_\edge + \tau_\edge \Rdisc(\und, \rnd, \vec f_\edge, \vec g_\edge) }_\node - \vec g_\node\big) \cdot \psidisc_\node\big],
 \end{split}
\end{align*}
where we again set $\und = \rnd = 0$ for all $\node \in \setNode\setminus \setNodeDir$.
After having computed the solution~$(\lambdadisc,\phidisc)$ to the condensed HDG formulation~\cref{EQ:hybrid_problem_disc}, the HDG solution can be retrieved as follows: The hybrid HDG unknowns are obtained by $\undisc = \lambdadisc_\node + \und$ and $\rndisc = \phidisc_\node + \rnd$ for all $\node \in \setNode$, and the primal and dual HDG unknowns $\uedisc, \rndisc,\ndisc,\mdisc$ can be obtained for all $\edge \in \setEdge$ similar to \cref{EQ:retrievsolution}, but with the local mappings \cref{EQ:disclocmapping}.

\begin{remark}[Size of the global system of equations]
 Importantly, the global system of equations corresponding to \cref{EQ:hybrid_problem_disc} does not change its size if the polynomial degree $p$ is increased. This is because the hybrid unknowns live in a zero-dimensional domain, and a single parameter can characterize all polynomials in such a domain.
\end{remark}

The following lemma also shows the symmetry and positive definiteness of the bilinear form $\bar A$, which proves the well-posedness of the condensed HDG  formulation.
\begin{lemma}[Properties of condensed HDG problem]\label{LEM:disc_spd_form}
 The bilinear form $\bar A$ can be equivalently rewritten as the following symmetric expression
 \begin{align*}
  &\bar A((\lambdadisc, \phidisc),(\mudisc, \psidisc))\\ 
  &\qquad = \sum_{\edge \in \setEdge}\big[ (C^{-1}_{\vec n}  \bar{\vec N}_\edge(\lambdadisc_\node, \phidisc_\node),  \bar{\vec N}_\edge(\mudisc_\node, \psidisc_\node))_\edge + (C^{-1}_{\vec m}  \Mdisc(\lambdadisc_\node, \phidisc_\node),  \Mdisc(\mudisc_\node, \psidisc_\node))_\edge\big]\\
  &\qquad \qquad \quad+ \sum_{\node \in \setNode\setminus \setNodeDir} \big[\big(\jump{  \tau_\edge (\bar{\vec U}_\edge(\lambdadisc_\node, \phidisc_\node) - \lambdadisc_\node) \cdot (\bar{\vec U}_\edge(\mudisc_\node, \psidisc_\node) - \mudisc_\node) }_\node\\[-2ex]
  &\hspace{3.75cm}\quad
    + \jump{  \tau_\edge (\Rdisc(\lambdadisc_\node, \phidisc_\node) - \phidisc_\node) \cdot (\Rdisc(\mudisc_\node, \psidisc_\node) - \psidisc_\node) }_\node\big)\big].
 \end{align*}
 Therefore, for any $\tau_\edge>0$, the bilinear form $\bar A$ is symmetric and positive definite. Consequently, the condensed formulation \cref{EQ:hybrid_problem_disc} of the HDG method is well-posed.
\end{lemma}

\begin{proof}
 The proof of this result uses very similar arguments as the proof of \cref{LEM:characterize_bil}, but applied to~\cref{EQ:discrete_solver} instead of~\cref{EQ:local_solver}. Therefore, it is omitted for the sake of brevity.
\end{proof}

\section{Error analysis}\label{SEC:error}
%
In this section, we perform an a priori error analysis of the proposed HDG discretization for the Timoshenko beam network problem. The error analysis is inspired by \cite{CockburnGS10}, where a projection-based error analysis of the HDG method for an elliptic model problem is performed. The projection is tailored to the HDG structure, allowing for a relatively simple and concise error analysis. In the present setting of Timoshenko beam network problems, the edgewise defined projection is given by the mapping $\Pi \colon (H^1(\edge))^3 \times (H^1(\edge))^3 \to (\mathds P_p(\edge))^3 \times (\mathds P_p(\edge))^3$, whose components $\Pi(\vec u_\edge, \vec n_\edge) = (\Pi_1(\vec u_\edge, \vec n_\edge), \Pi_2(\vec u_\edge, \vec n_\edge))$ are for all $\edge \in \setEdge$ defined by the following conditions:
\begin{subequations}\label{EQ:projection}
 \begin{align}
  (\Pi_1(\vec u_\edge, \vec n_\edge), \bar{\vec v} )_\edge & = (\vec u_\edge, \bar{\vec v})_\edge & \text{ for all } \bar{\vec v} \in \mathds P_{p-1}(\edge),\label{EQ:projection_1} \\
  (\Pi_2(\vec u_\edge, \vec n_\edge), \bar{\vec v} )_\edge & = (\vec n_\edge, \bar{\vec v})_\edge & \text{ for all } \bar{\vec v} \in \mathds P_{p-1}(\edge), \label{EQ:projection_2}\\
  [\Pi_2(\vec u_\edge, \vec n_\edge)] \nu + \tau_\edge [\Pi_1(\vec u_\edge, \vec n_\edge)] & = \vec n_\edge \nu + \tau_\edge \vec u_\edge & \text{ for all } \node \sim \edge,\label{EQ:projection_3}
 \end{align}
\end{subequations}
where we recall that we write $\node \sim \edge$ if the node $\node$ is an endpoint of the edge $\edge$. The well-posedness of this projection can be concluded from~\cite[Thm.~2.1]{CockburnGS10}, where the result was proved in a more general multi-dimensional setting.

To simplify the notation of the following error analysis, we define  the errors
\begin{equation}\label{eq:defdiscandprojerr}
	\begin{aligned}
		 \theta^{\vec u}_\edge &= \vec u_\edge - \Pi_1(\vec u_\edge, \vec n_\edge),\qquad  & \epsilon^{\vec u}_\edge &= \Pi_1(\vec u_\edge, \vec n_\edge) - \bar{\vec u}_\edge,\qquad  & \theta^{\vec u}_\edge  + \epsilon^{\vec u}_\edge & = \vec u_\edge - \bar{\vec u}_\edge\\
		\theta^{\vec n}_\edge &= \vec n_\edge - \Pi_2(\vec u_\edge, \vec n_\edge), & \epsilon^{\vec n}_\edge &= \Pi_2(\vec u_\edge, \vec n_\edge) - \bar{\vec n}_\edge, & \theta^{\vec n}_\edge + \epsilon^{\vec n}_\edge & = \vec n_\edge - \bar{\vec n}_\edge,
	\end{aligned}
\end{equation}
and similar definitions are used for the variables $\re$ and $\vec m_\edge$.

The following lemma states projection error estimates for the operator $\Pi$, which are taken from \cite[Thm.~2.1]{CockburnGS10}. We denote the $H^k(\edge)$-seminorm for any $k \in \mathbb N_0$ by $|\cdot|_{\edge,k}$. 

\begin{lemma}[Projection error estimates]\label{LEM:proj_bound}
	Given $p_{\vec u}, p_{\vec n} \in [0,p]$, there holds for any edge $\edge \in \setEdge$ and all $\vec u_\edge \in H^{p_{\vec u}+1}(\edge)$, $\vec n_\edge \in H^{p_{\vec n}+1}(\edge)$ that 
 \begin{align*}
  \| \theta^{\vec u}_\edge \|_\edge & \lesssim \edgelength^{p_{\vec u} +1} | \vec u_\edge |_{\edge,p_{\vec u}+1} + \tfrac1{\tau_\edge} \edgelength^{p_{\vec n}+1} | \vec n_\edge |_{\edge,p_{\vec n}+1}
   \\
  \| \theta^{\vec n}_\edge \|_\edge & \lesssim \tau_\edge \edgelength^{p_{\vec u}+1} | \vec u_\edge |_{\edge,p_{\vec u}+1} + \edgelength^{p_{\vec n}+1} | \vec n_\edge |_{\edge,p_{\vec n}+1},
 \end{align*}
 and analogous estimates hold for $\theta^{\vec r}_\edge$ and $\theta^{\vec m}_\edge$.
\end{lemma} 

Recall that the primal unknowns are edgewise $H^1$-regular and satisfy continuity conditions at the nodes. Therefore, for any edge $\edge \in \setEdge$ and any node $\node\sim \edge$, $\vec u_\edge(\node) = \vec u_\node$ and $\vec r_\edge(\node) = \vec r_\node$ holds. This motivates omitting the primal unknowns' subscripts by simply writing $\vec u$ and $\vec r$. Similarly, we may omit the subscript of $\vec n_\edge$ and $\vec m_\edge$. The following lemma states error equations that will be frequently used in the subsequent error analysis. 

\begin{lemma}[Error equations]
	
 For any edge $\edge \in \setEdge$ and all test functions $\bar{\vec p}, \bar{\vec q}, \bar{\vec v}, \bar{\vec w} \in (\mathbb P_p(\edge))^3$, there hold the identities
 \begin{subequations}
 \begin{align}
  (\epsilon^{\vec u}_\edge, \partial_x \bar{\vec p})_\edge & =
  (C^{-1}_{\vec n} (\vec n - \bar{\vec n}_\edge), \bar{\vec p})_\edge+ (\vec i_\edge \times (\vec r - \redisc), \bar{\vec p})_\edge + \langle\vec u - \bar{\vec u}_\node,\bar{\vec p}\nu_\edge\rangle_\edge,  \label{EQ:eu_ident}\\
  (\epsilon^{\vec r}_\edge, \partial_x \bar{\vec q})_\edge & = (C^{-1}_{\vec m}(\vec m - \bar{\vec m}_\edge), \bar{\vec q})_\edge + \langle\vec r - \rndisc,\bar{\vec q} \nu_\edge\rangle_\edge, \label{EQ:ephi_ident}\\
  (\partial_x \epsilon^{\vec n}_\edge, \bar{\vec v})_\edge & = \langle\tau_\edge (\bar{\vec u}_\edge - \bar{\vec u}_\node) - \theta^{\vec n}_\edge  \nu_\edge,\bar{\vec v}\rangle_\edge, \label{EQ:en_ident}\\
  (\partial_x \epsilon^{\vec m}_\edge, \bar{\vec w})_\edge & 
  = \langle \tau_\edge (\redisc - \rndisc) - \theta^{\vec m}_\edge \nu_\edge,\bar{\vec w}\rangle_\edge - (\vec i_\edge \times (\vec n - \bar{\vec n}_\edge), \bar{\vec w})_\edge.\label{EQ:em_ident}
 \end{align}
 Furthermore, for any nodal functions $\vec v_\node, \vec w_\node$ with support in $\node \in \setNode \setminus \setNodeDir$ we have that
 \begin{equation}
  0 =  \jump{ \epsilon_\edge^{\vec n} \nu_\edge + \tau_\edge \epsilon_\edge^{\vec u} - \tau_\edge (\vec u-\bar{\vec u}_\node) }_\node \vec v_\node + \jump{ \epsilon_\edge^{\vec m}\nu_\edge + \tau_\edge \epsilon_\edge^{\vec r} - \tau_\edge (\vec r-\rndisc) }_\node \vec w_\node .\label{EQ:error_fluxes}
 \end{equation}
 \end{subequations}
\end{lemma}
\begin{proof}
 Error equations \cref{EQ:eu_ident,EQ:ephi_ident} can be derived from \cref{EQ:local_solver,EQ:discrete_solver}, using property \cref{EQ:projection_1} of the projection $\Pi$. To obtain \cref{EQ:en_ident,EQ:em_ident}, an integration by parts is necessary before property \cref{EQ:projection_2} of $\Pi$ can be used. The error equation \cref{EQ:error_fluxes} can be derived using property \cref{EQ:projection_3} of $\Pi$ and the balance conditions \cref{EQ:timo_balance,EQ:timo_balance_disc}. 
\end{proof}

The following two lemmas provide preliminary results, combined in \cref{COR:errorest} to a full error estimate of the proposed HDG discretization. The first lemma provides an estimate for the dual unknowns.

\begin{lemma}[Estimate for dual unknowns]\label{LEM:duals_conv}
 It holds that
  \begin{equation}
  	\label{EQ:estenem}
 	\sum_{\edge \in \setEdge}  \| \epsilon^{\vec n}_\edge\|^2_\edge 
 	\lesssim
 	\sum_{\edge \in \setEdge} \left[ \| \theta^{\vec n}_\edge\|^2_\edge +  \| \theta^{\vec r}_\edge \|^2_\edge\right], \qquad 
 	\sum_{\edge \in \setEdge}  \|\epsilon^{\vec m}_\edge\|^2_\edge 
 	\lesssim  \sum_{\edge \in \setEdge} \left[  \| \theta^{\vec m}_\edge\|^2_\edge+\|\theta^{\vec n}_\edge\|_\edge\|\epsilon^{\vec r}_\edge\|_\edge \right].
 \end{equation}
\end{lemma}
\begin{proof}
 We begin this proof with the obvious identity
 \begin{equation}\label{EQ:obviousidentity}
  - \langle \epsilon^{\vec n}_\edge\nu_\edge + \tau_\edge \epsilon^{\vec u}_\edge - \tau_\edge(\vec u-\undisc), \vec u -\undisc \rangle_\edge= - \langle \epsilon^{\vec n}_\edge\nu, \vec u -\undisc \rangle_\edge + \langle \tau_\edge(\vec u-\undisc)-\tau_\edge \epsilon^{\vec u}_\edge, \vec u -\undisc \rangle_\edge.
 \end{equation}
 To derive an equivalent representation for the left-hand side of \cref{EQ:obviousidentity}, we test~\cref{EQ:eu_ident} with the test function~$\epsilon^{\vec n}_\edge$ and rearrange the terms which yields that
 \begin{align}\label{EQ:tested}
 - \langle\vec u-\undisc,\epsilon^{\vec n}_\edge\nu_\edge\rangle_\edge=   (C_{\vec n}^{-1}(\vec n -\ndisc),\epsilon^{\vec n}_\edge)_\edge + (\vec i_\edge \times (\vec r - \redisc),\epsilon^{\vec n}_\edge)_\edge -(\epsilon^{\vec u}_\edge ,\partial_x \epsilon^{\vec n}_\edge)_\edge.
 \end{align}
 The last term on the right-hand side of \cref{EQ:tested} can be further rewritten using the identity 
 \begin{align}\label{EQ:tested2}
  (\partial_x \epsilon^{\vec n}_\edge,\epsilon^{\vec u}_\edge)_\edge =\langle\tau_\edge (\bar{\vec u}_\edge - \bar{\vec u}_\node) - \theta^{\vec n}_\edge \nu_\edge,\epsilon^{\vec u}_\edge\rangle_\edge = \langle\tau_\edge(\vec u-\undisc)-\tau_\edge \epsilon^{\vec u}_\edge ,\epsilon^{\vec u}_\edge\rangle_\edge,
 \end{align}
 which can be derived by testing \cref{EQ:en_ident} with the test function $\epsilon^{\vec u}_\edge$ and using property \cref{EQ:projection_3} of the projection~$\Pi$. Inserting \cref{EQ:tested,EQ:tested2} into \cref{EQ:obviousidentity} then yields that
 \begin{equation}\label{EQ:identity_jumpu}
  \begin{aligned}
   - &  \langle \epsilon^{\vec n}_\edge\nu_\edge + \tau_\edge \epsilon^{\vec u}_\edge - \tau_\edge(\vec u-\undisc), \vec u -\undisc \rangle_\edge\\
   &= (C_{\vec n}^{-1}(\vec n -\ndisc),\epsilon^{\vec n}_\edge)_\edge + (\vec i_\edge \times (\vec r - \redisc),\epsilon^{\vec n}_\edge)_\edge +\tau_\edge \langle \epsilon^{\vec u}_\edge -(\vec u -\undisc), \epsilon^{\vec u}_\edge -(\vec u -\undisc) \rangle_\edge,
  \end{aligned}
 \end{equation}
 which is the desired equivalent representation.

 Summing the latter identity over all edges $\edge \in \setEdge$ and using \cref{EQ:error_fluxes} for the nodal functions ${{\vec v}_\node = \vec u - \undisc}$, we obtain that
 \begin{equation}\label{EQ:summed1}
  0 = \sum_{\edge \in \setEdge} \left[ (C_{\vec n}^{-1}(\vec n -\ndisc),\epsilon^{\vec n}_\edge)_\edge + (\vec i_\edge \times (\vec r - \redisc),\epsilon^{\vec n}_\edge)_\edge \right] + \sum_{\node \in \setNode} \tau_\edge \jump{(\epsilon^{\vec u}_\edge -(\vec u -\undisc))^2}_\node,
 \end{equation}
where we interpret the square of a vector as a dot product with itself. Algebraic manipulations and the non-negativity of the last term in \cref{EQ:summed1} then yield the estimate
 \begin{equation*}
	\sum_{\edge \in \setEdge} (C_{\vec n}^{-1}\epsilon^{\vec n}_\edge,\epsilon^{\vec n}_\edge)_\edge  \leq \sum_{\edge \in \setEdge} \left[ -(C_{\vec n}^{-1}\theta^{\vec n}_\edge,\epsilon^{\vec n}_\edge)_\edge  - (\vec i_\edge \times \theta^{\vec r}_\edge,\epsilon^{\vec n}_\edge)_\edge \right].
\end{equation*}
The first inequality in \cref{EQ:estenem} can then be concluded from H\"older's and Young's inequalities and the uniform bounds of the coefficients $C_{\vec n}$ and $C_{\vec m}$, cf.~\cref{eq:spectralbounds}. 

Similar considerations can be applied to the rotations and moments, which yields that
 \begin{align}\label{EQ:summed2}
  0 = \sum_{\edge \in \setEdge} \left[ (C_{\vec m}^{-1}(\vec m -\mdisc),\epsilon^{\vec m}_\edge)_\edge + (\vec i_\edge \times (\vec n - \ndisc),\epsilon^{\vec r}_\edge)_\edge \right] + \sum_{\node \in \setNode} \tau_\edge \jump{ (\epsilon^{\vec r}_\edge -(\vec r -\rndisc))^2 }_\node,
 \end{align}
and the second inequality in \cref{EQ:estenem} can again be obtained by using the non-negativity of the last term in \cref{EQ:summed2}, as well as H\"older's and Young's inequalities and bounds~\cref{eq:spectralbounds}. 
\end{proof}

The second lemma provides an estimate for the primal unknowns.
\begin{lemma}[Estimate for primal unknowns]\label{LEM:primalest}
	 It holds that
 \begin{align*}
 	\sum_{\edge \in \setEdge} \big[\| \epsilon^{\vec u}_\edge  \|^2_\edge + \| \epsilon^{\vec r}_\edge  \|^2_\edge\big] \lesssim 
 \sum_{\edge \in \setEdge} (\edgelength+\tau_\edge \edgelength^2)^2&\left[ \| \theta^{\vec n}_\edge \|^2_\edge + \| \theta^{\vec m}_\edge \|^2_\edge + \| \theta^{\vec r}_\edge \|^2_\edge \right]\\
+  \sum_{\edge \in \setEdge}\big(\tfrac{\edgelength}{\tau_\edge} \big)^2 &\left[\|\theta^{\vec n}_\edge\|_\edge^2 + \|\theta^{\vec r}_\edge\|_\edge^2\right].
 \end{align*}
\end{lemma}

\begin{proof}
This proof is based on an Aubin--Nietsche-type argument. For this, we consider an auxiliary problem whose data and solution we indicate by a dagger. The auxiliary problem is the same as problem~\cref{EQ:timo_network}, but with homogeneous concentrated forces and moments and homogeneous Dirichlet data, i.e., $\vec f_\node^\dagger=\vec g_\node^\dagger = 0$ for all~$\node \in \setNode$ and $\vec u_\node^\dagger = \vec r^\dagger_\node = 0$ for all~$\node \in \setNodeDir$. The following elliptic regularity estimate holds for the auxiliary problem:
\begin{equation}\label{EQ:ellregauxiliaryproblem}
	\sum_{\edge \in \setEdge} \left[ | \vec u^\dagger_\edge |^2_{2,\edge} + | \re^\dagger |^2_{2,\edge} + | \vec n^\dagger_\edge |^2_{1,\edge} + | \vec m^\dagger_\edge |^2_{1,\edge} \right] \lesssim \sum_{\edge \in \setEdge} \left[ \| \vec f^\dagger_\edge \|^2_\edge + \| \vec g^\dagger_\edge \|^2_\edge \right].
\end{equation}
Note that since the auxiliary problem satisfies the same continuity properties as the original problem, we may also write $\vec u^\dagger$ for~$\vec u_\edge^\dagger$ and~$\vec u_\node^\dagger$, and similarly for $\vec r^\dagger$, $\vec n^\dagger$, and $\vec m^\dagger$.

We begin the proof by deriving an equivalent expression for $(\epsilon_\edge^{\vec u},\vec f_\edge^\dagger)_\edge$. By the definition of $\vec n^\dagger$, cf.~\cref{EQ:timo_force}, integration by parts, and properties~\cref{EQ:projection_2,EQ:projection_3}, we get that
 \begin{align}
  (\epsilon^{\vec u}_\edge,\vec f^\dagger_\edge)_\edge & = (\epsilon^{\vec u}, \partial_x \Pi_2(\vec u^\dagger, \vec n^\dagger))_\edge + \langle \epsilon^{\vec u}_\edge, \tau_\edge (\Pi_1(\vec u^\dagger, \vec n^\dagger) - \vec u^\dagger)\rangle_\edge \notag\\
 \begin{split}\label{EQ:T2}
  & = (C^{-1}_{\vec n}(\vec n - \bar{\vec n}_\edge), \Pi_2(\vec u^\dagger, \vec n^\dagger))_\edge + (\vec i_\edge \times (\vec r - \bar{\vec r}_\edge), \Pi_2(\vec u^\dagger, \vec n^\dagger))_\edge \\
  & \qquad\quad  + \langle \vec u - \bar{\vec u}_\node, \Pi_2(\vec u^\dagger, \vec n^\dagger) \nu_\edge\rangle_\edge + \langle \epsilon^{\vec u}_\edge , \tau_\edge (\Pi_1(\vec u^\dagger, \vec n^\dagger) - \vec u^\dagger)\rangle_\edge,
 \end{split}
 \end{align}
 where we used~\cref{EQ:eu_ident} in the last equality. 

In the following, we will sum equation~\cref{EQ:T2} over all edges $\edge \in\setEdge$. For the last two terms in~\cref{EQ:T2}, this results in
 \begin{align}
 &\sum_{\edge \in \setEdge}\big[\langle \vec u - \bar{\vec u}_\node, \Pi_2(\vec u^\dagger, \vec n^\dagger) \nu_\edge\rangle_\edge + \langle \epsilon^{\vec u}_\edge , \tau_\edge (\Pi_1(\vec u^\dagger, \vec n^\dagger) - \vec u^\dagger)\rangle_\edge\big]\notag\\
		&\qquad =\sum_{\node \in \setNode} \jump{(\vec u - \bar{\vec u}_\node) \cdot (\Pi_2(\vec u^\dagger, \vec n^\dagger) - \vec n^\dagger) \nu_\edge } + \sum_{\node \in \setNode} \jump{ \epsilon^{\vec u}_\edge \cdot \tau_\edge [\Pi_1(\vec u^\dagger, \vec n^\dagger) - \vec u^\dagger]}\notag \\
		& \qquad= \sum_{\node \in \setNode} \jump{ [\epsilon^{\vec u}_\edge - (\vec u - \bar{\vec u}_\node)] \cdot \tau_\edge [\Pi_1(\vec u^\dagger, \vec n^\dagger) - \vec u^\dagger]} 
		=  \sum_{\edge \in \setEdge} (\vec \epsilon^{\vec n}_\edge, \partial_x \vec u^\dagger)_\edge.\label{EQ:identitysum}
\end{align}
Note that the first equality follows after rearranging the summands and inserting the term~$\jump{(\vec u-\undisc)\cdot \vec n^\dagger\nu_\edge}_\node$ which equals zero since $\vec u - \undisc=0$ for all $\node \in \setNodeDir$ and $\jump{\vec n^\dagger\nu_\edge}_\node = 0$ for all $\node\in \setNode\setminus \setNodeDir$. The second equality is a direct consequence of \cref{EQ:projection_3}.  To prove the last equality, we use the following two identities
 \begin{align}
 	 	\sum_{\node \in \setNode} \jump{ [\epsilon^{\vec u}_\edge - (\vec u - \bar{\vec u}_\node)] \cdot \tau_\edge \Pi_1(\vec u^\dagger, \vec n^\dagger)}_\node &= -\sum_{\node \in \setNode}\jump{[\tau_\edge (\uedisc-\undisc) - \theta^{\vec n}_\edge \nu_\edge] \cdot \Pi_1(\vec u^\dagger,\vec n^\dagger)}_\node,\label{EQ:ident1}\\
 	\sum_{\node \in \setNode} \jump{ [\epsilon^{\vec u}_\edge - (\vec u - \bar{\vec u}_\node)] \cdot \tau_\edge \vec u^\dagger}_\node &=-\sum_{\node \in \setNode} \jump{ \epsilon^{\vec \node}_\edge\nu_\edge \cdot \vec u^\dagger}_\node,\label{EQ:ident2}
 \end{align}
 as well as \cref{EQ:en_ident} and integrate by parts. Identity \cref{EQ:ident1} can be concluded from~\cref{EQ:projection_3}, and  \cref{EQ:ident2} follows from \cref{EQ:projection_3} and the balance condition $\jump{\ndisc \nu_\edge + \tau_\edge (\uedisc-\undisc)}_\node = \jump{\vec n_\edge \nu_\edge}_\node = \vec f_\node$. 
 
 Denoting by $\pi$ the edgewise $L^2$-orthogonal projection onto polynomials of degree at most $p-1$, we finally obtain for the sum of equation~\cref{EQ:T2} over all edges $\edge \in \setEdge$ that
  \begin{align*}
  \sum_{\edge \in \setEdge} (\epsilon^{\vec u}_\edge, \vec f^\dagger_\edge)_\edge 
  & = \sum_{\edge \in \setEdge} (C^{-1}_{\vec n}(\vec n - \bar{\vec n}_\edge),  \Pi_2(\vec u^\dagger, \vec n^\dagger) - \vec n^\dagger)_\edge + \sum_{\edge \in \setEdge}(\vec i_\edge \times (\vec r - \bar{\vec r}_\edge), \Pi_2(\vec u^\dagger, \vec n^\dagger))_\edge \\
  & \qquad - \sum_{\edge \in \setEdge} (\vec n - \bar{\vec n}_\edge, \vec i_\edge \times \vec r^\dagger)_\edge +  \sum_{\edge \in \setEdge} (\Pi_2(\vec u, \vec n)-\vec n, \partial_x \vec u^\dagger - \pi \partial_x \vec u^\dagger)_\edge,
 \end{align*}
 where we used \cref{EQ:identitysum} as well as the identity
\begin{equation*}
	(C^{-1}_{\vec n} (\vec n - \bar{\vec n}_\edge), \vec n^\dagger)_\edge = -(\vec n - \bar{\vec n}_\edge,\partial_x\vec u^\dagger+\vec i_\edge \times \vec r^\dagger)_\edge,
\end{equation*}
which follows by the definition of $\vec n^\dagger$, cf.~\cref{EQ:timo_network}. Note that the term $\pi \partial_x \vec u^\dagger$ does not contribute anything since it is a polynomial of degree at most $p-1$. 

 Similar considerations can be applied to derive an equivalent representation for the sum of~$(\epsilon_\edge^{\vec r},\vec g_\edge^\dagger)$ over all edges $\edge \in \setEdge$, which results in
 \begin{align*}
  \sum_{\edge \in \setEdge} (\epsilon^{\vec r}_\edge, \vec g^\dagger_\edge)_\edge & = \sum_{\edge\in \setEdge} (C^{-1}_{\vec m} (\vec m - \bar{\vec m}_\edge),  \Pi_2(\vec r^\dagger, \vec m^\dagger) - \vec m^\dagger)_\edge + \sum_{\edge\in \setEdge}( \epsilon^{\vec r}_\edge , \vec i_\edge \times \vec n^\dagger)_\edge \\
  & \qquad - \sum_{\edge \in \setEdge}(\vec i_\edge \times (\vec n - \bar{\vec n}_\edge), \Pi_1(\vec r^\dagger,\vec m^\dagger))_\edge + \sum_{\edge\in \setEdge} (\Pi_2(\vec r, \vec m)-\vec m, \partial_x \vec r^\dagger - \pi \partial_x \vec r^\dagger)_\edge.
 \end{align*}

Combining the above representations, we obtain, after some algebraic manipulation,~that
 \begin{align}
  \sum_{\edge \in \setEdge}& \big[ (\epsilon^{\vec u}_\edge, \vec f^\dagger_\edge)_\edge  +(\epsilon^{\vec r}_\edge, \vec g^\dagger_\edge)_\edge\big]\notag\\
     =& -\sum_{\edge \in \setEdge} \big[(C^{-1}_{\vec n}( \epsilon^{\vec n}_\edge + \theta^{\vec n}_\edge), \theta^{\vec n^\dagger}_\edge)_\edge + (C^{-1}_{\vec m}( \epsilon^{\vec m}_\edge + \theta^{\vec m}_\edge), \theta^{\vec m^\dagger}_\edge)_\edge\big] \notag\\
  &- \sum_{\edge \in \setEdge}\big[ (\theta^{\vec n}_\edge, \partial_x \vec u^\dagger - \pi \partial_x \vec u^\dagger)_\edge +(\theta^{\vec m}_\edge, \partial_x \vec r^\dagger - \pi\partial_x \vec r^\dagger_\edge)_\edge\big] + \sum_{\edge \in \setEdge} ( \vec i_\edge \times \theta^{\vec r}_\edge, \vec n^\dagger - \pi \vec n^\dagger)\notag\\
    &-\sum_{\edge \in \setEdge}\big[ ( \vec i_\edge \times (\theta^{\vec r}_\edge+\epsilon_\edge^{\vec r}), \theta^{\vec n^\dagger}_\edge)-(\vec i_\edge \times (\theta^{\vec n}_\edge+\epsilon_\edge^{\vec n}),\theta^{\vec r^\dagger}_\edge)_\edge\big],\label{EQ:eqrep}
 \end{align}
 where the term $\pi \vec n^\dagger$ could be inserted without contribution since it is a polynomial of degree at most $p-1$.
 
 Using the identity of norms
 \begin{equation*}
 	\bigg[\sum_{\edge \in \setEdge} \big[\| \epsilon^{\vec u}_\edge  \|^2_\edge + \| \epsilon^{\vec r}_\edge  \|^2_\edge\big]  \bigg]^{1/2}=  \sup_{\sum_{\edge \in \setEdge}\|\vec f^\dagger_\edge\|_\edge^2 + \|\vec g^\dagger_\edge\|_\edge^2 = 1}\;\; \sum_{\edge \in \setEdge} \big[(\epsilon^{\vec u}_\edge, \vec f^\dagger_\edge)_\edge +  (\epsilon^{\vec r}_\edge, \vec g^\dagger_\edge)_\edge\big],
 \end{equation*}
which is a consequence of the Riesz representation theorem, as well as~\cref{EQ:eqrep}, \cref{LEM:duals_conv,LEM:proj_bound}, classical approximation properties of the $L^2$-projection $\pi$, and the elliptic regularity result~\cref{EQ:ellregauxiliaryproblem}, we obtain with Young's and Hölder's inequalities that
 \begin{align*}
 	\sum_{\edge \in \setEdge} \big[\| \epsilon^{\vec u}_\edge  \|^2_\edge + \| \epsilon^{\vec r}_\edge  \|^2_\edge\big]
 	 \lesssim  \sum_{\edge \in \setEdge} (\edgelength+\tau_\edge \edgelength^2)^2&\left[ \| \theta^{\vec n}_\edge \|^2_\edge + \| \theta^{\vec m}_\edge \|^2_\edge + \| \theta^{\vec r}_\edge \|^2_\edge \right] +   \sum_{\edge \in \setEdge}(\edgelength+\tau_\edge \edgelength^2)^2 \| \epsilon^{\vec r}_\edge  \|^2_\edge \\
 	 + \sum_{\edge \in \setEdge}\big(\tfrac{\edgelength}{\tau_\edge} \big)^2 &\left[\|\theta^{\vec n}_\edge\|_\edge^2 + \|\theta^{\vec r}_\edge\|_\edge^2\right].
 \end{align*}
For sufficiently small edge lengths $\edgelength$, the assertion can be concluded by absorbing the term involving $\epsilon_\edge^{\vec r}$ in the left-hand side, where we recall that $\tau_\edge \edgelength \lesssim 1$. 
\end{proof}

We are now ready to present the desired convergence result for the proposed HDG discretization. It provides estimates for the $L^2$-errors of the primal and dual variables in terms of the maximum edge length, denoted by 
\begin{equation*}
	\maxedgelength \coloneqq \max_{\edge\in \setEdge} \edgelength.
\end{equation*}
The rates we prove agree with those found in classical HDG theory, cf.~\cite{CockburnGS10}.

\begin{theorem}[Convergence of HDG method]\label{COR:errorest}
	Suppose that for any edge $\edge \in \setEdge$ it holds that $\vec u_\edge, \vec r_\edge, \vec n_\edge, \vec m_\edge \in H^{p+1}(\edge)$. Then, given a stabilization parameter which scales like~$\edgelength^s$ for some $s \in \{-1,0,1\}$, the HDG approximation converges to the solution of the Timoshenko beam network model with the error estimates  
	\begin{alignat}{2}
		\bigg[\sum_{\edge \in \setEdge}\left[\|\vec u_\edge  - \bar{\vec u}_\edge\|_\edge^2 + \|\vec r_\edge - \bar{\vec r}_\edge\|_\edge^2\right]\bigg]^{1/2} &\lesssim \maxedgelength^{p+1-s^+},\label{EQ:L2primal}\\
	\bigg[\sum_{\edge \in \setEdge}\left[\|\vec n_\edge - \bar{\vec n}_\edge\|_\edge^2 + \|\vec m_\edge - \bar{\vec m}_\edge\|_\edge^2\right]\bigg]^{1/2} &\lesssim \maxedgelength^{p+1-|s|},\label{EQ:L2dual}
\end{alignat}
where we denote $s^+ \coloneqq \max(s,0)$.
\end{theorem}
\begin{proof}
	To prove the above error estimates, we use that the error can be written as the sum of a projection and a discrete error, cf. \cref{eq:defdiscandprojerr}. For estimating the projection error we use \cref{LEM:proj_bound}, while for the discrete error \cref{LEM:primalest,LEM:duals_conv} can be used. For example, for the choice $\tau_\edge \sim 1$, this gives a convergence of order $p+1$ for both the projection error and the discrete errors, and thus an overall convergence of order $p+1$, cf. \cref{EQ:L2primal,EQ:L2dual}.
\end{proof}
Note that using the convergence result for the dual variable, it seems possible to follow classical HDG postprocessing techniques to boost the convergence of the primal variables to order $p+2$.

\section{Domain decomposition preconditioner}\label{SEC:precond}

A practical implementation of the proposed HDG method requires solving the linear system of equations corresponding to~\cref{EQ:hybrid_problem_disc}. The linear systems of equations that typically arise when simulating, for example, the structural properties of paper-based materials, cf.~\cref{fig:paper}, can easily become intractable for direct solvers due to their large size. In addition, they are typically poorly conditioned, requiring appropriate preconditioners. This section introduces a two-level overlapping additive Schwarz preconditioner similar to~\cite{GoHeMa22}, which we will use within a preconditioned conjugate gradient method. 

The domain decomposition and coarse space used by the proposed preconditioner are constructed using an artificial (coarse) mesh $\mathcal{T}_H$ of a bounding domain $\Omega\subset \mathbb{R}^3$ of the spatial network. 
\begin{figure}
 \includegraphics[height=.375\textwidth]{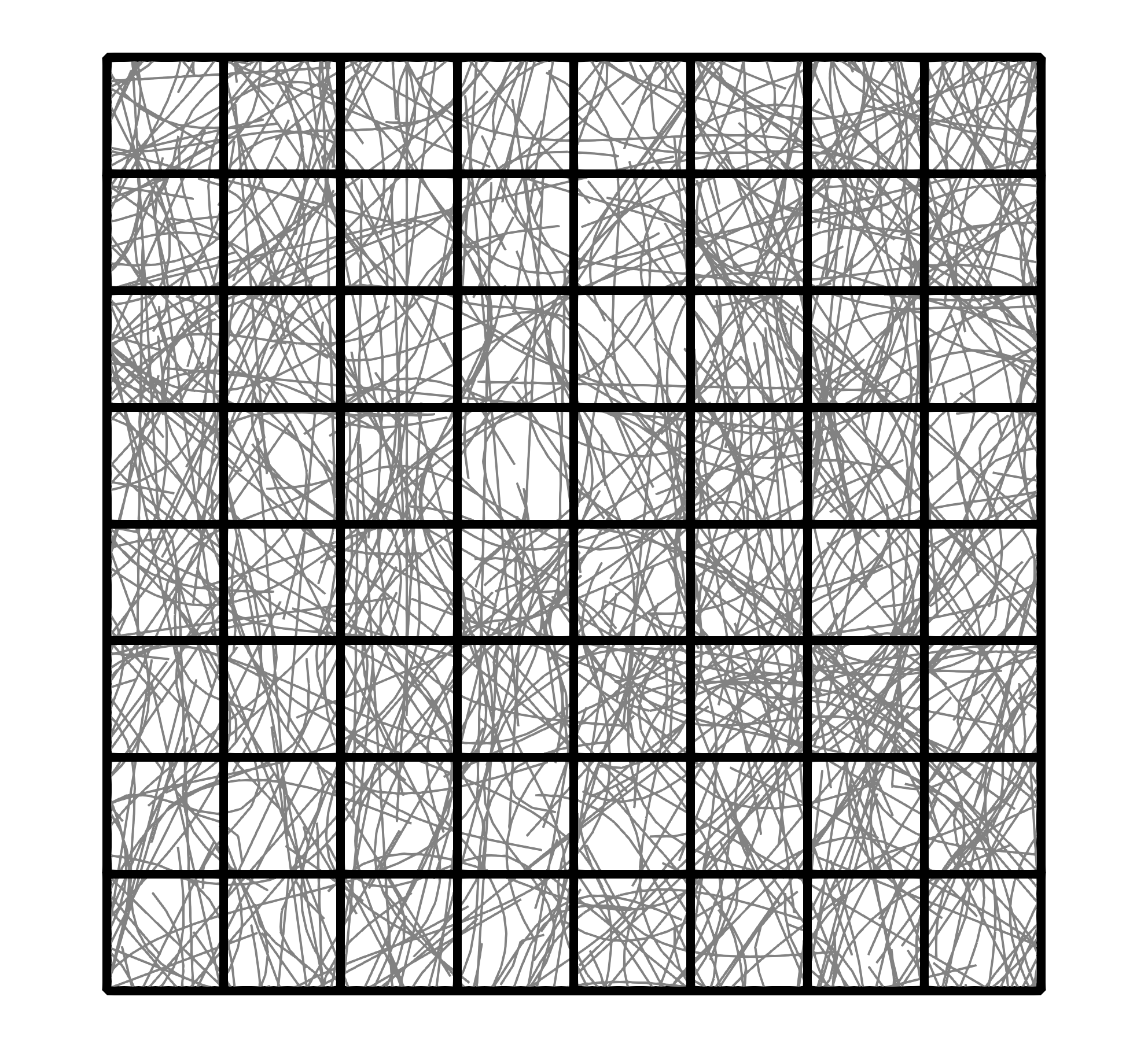}\hspace{.5cm}
 \includegraphics[height=.375\textwidth]{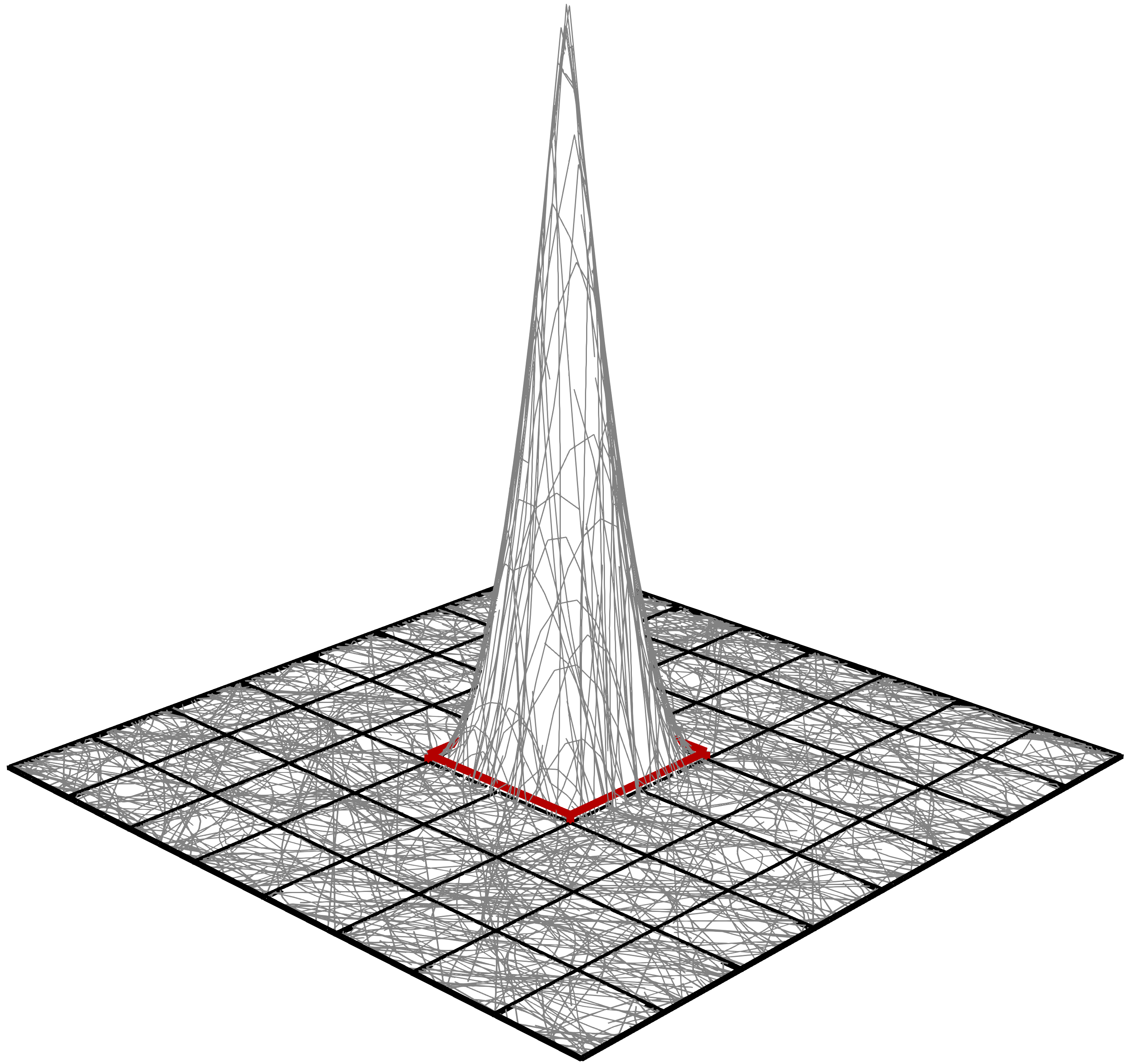}
 \caption{An artificial mesh $\mathcal{T}_H$ over a network (left) and a basis function~$\varphi_i$ with the boundary of its support marked in red (right).}\label{fig:artificialgrid}
\end{figure}
For simplicity, we will assume that $\Omega$ is a box equipped with a uniform Cartesian mesh.  
Note that for thin materials such as cardboard, the corresponding artificial mesh may contain much fewer elements in one spatial dimension than in the others.
With respect to the artificial mesh, we then introduce the set of trilinear basis functions $\{\varphi_i\}_{i=1}^m$, where $m$ is the number of nodes of the artificial mesh, and corresponding supports~$U_i\coloneqq \operatorname{supp}(\varphi_i)$. An illustration of an artificial mesh and corresponding basis functions in a two-dimensional setting can be found in \cref{fig:artificialgrid}.  The space of continuous piecewise trilinear functions with respect to $\mathcal{T}_H$, satisfying Dirichlet boundary conditions on boundary segments where the network nodes are fixed, is denoted by $V_H$. Note that the domain of the functions in~$V_H$ is considered to be the nodes of the spatial network. 

We employ a preconditioner based on the subspace decomposition
\begin{equation*}
	V_{\vec \lambda}=V_{\vec \lambda,0}+ V_{\vec \lambda,1}+\dots+ V_{\vec \lambda,m}
\end{equation*}
with the coarse space $V_{\vec \lambda,0}\coloneqq V_H\times V_H \times V_H$ and  local subspaces defined for $i=1,\dots,m$ by $V_{\vec \lambda,i}\coloneqq \{\vec v\in V_{\vec \lambda}\with \operatorname{supp}(\vec v)\subset U_i\}$. Given this decomposition, we introduce for any subspace a corresponding subspace projection operator $P_i\colon V_{\vec \lambda}\times V_{\vec \lambda}\rightarrow V_{\vec \lambda,i}\times V_{\vec \lambda,i}$ such~that
\begin{equation*}
	(\bar A P_i (\vec \lambda,\vec \phi),(\vec \mu,\vec \psi))=(\bar A (\vec \lambda,\vec \phi),(\vec \mu,\vec \psi)) 
\end{equation*}
holds for all $(\vec \mu,\vec \psi)\in V_{\vec \lambda,i}\times V_{\vec \lambda,i}$. Note that the existence and uniqueness of such an operator is a direct consequence of $\bar A$ being an inner product on $V_{\vec \lambda}\times V_{\vec \lambda}$. A preconditioned version of the operator $\bar A$ can then be defined as follows:
\begin{equation*}
	P\coloneqq P_0+P_1+\dots+P_m.
\end{equation*}
The preconditioner, denoted by $B$, is then given by the relation $P=B\bar A$. Note that the preconditioner is never explicitly formed. In practice, only the preconditioned operator must be computed, which requires the direct solution of a coarse global problem and~$m$ local problems, all of which can be solved independently. Because this approach involves direct solves, it is sometimes called semi-iterative. In the following, we will use this precondition for the conjugate gradient method.
  
The uniform convergence of the resulting preconditioned conjugate gradient method was proved in \cite[Thm.~4.3]{GoHeMa22} under the assumption that the considered linear system of equations is spectrally equivalent to the graph Laplacian, and under certain homogeneity, connectivity, and locality assumptions on the spatial network at coarse scales. 
The proof is inspired by classical Schwarz theory (see, e.g.,~\cite{Xu92,ToW05,KorY16}) and constructs a quasi-interpolation operator in the spatial network setting, whose approximation and stability properties could be proved using Poincare's and Friedrichs' inequalities on subgraphs. In practice, this preconditioner has demonstrated its ability to cope with the typically complex geometry of spatial networks and highly varying material properties; see, e.g.,~\cite{Grtz2024}. We emphasize that, in contrast, for example standard black-box preconditioners from the class of algebraic multigrid methods, cf.~\cite{XZ17,LB12}, may perform very poorly when applied to spatial network problems.  This is because they do not sufficiently consider the problem's geometry. 

To formulate the aforementioned spectral equivalence assumption required to prove the convergence of the preconditioned iteration, we introduce two bilinear forms acting on the space~$V_{\vec \lambda}$. The first bilinear form is a mass-type operator, while the second is a weighted graph Laplacian operator. They are for any functions  $\vec \lambda, \vec \mu \in V_{\vec \lambda}$ defined~by
\begin{align}
	\label{EQ:graphlaplacian}
	M(\vec \lambda,\vec \mu ) \coloneqq  \sum_{\node \in \setNode}\frac12\sum_{\edge \sim \node} \vec \lambda_\node \vec \mu_\node \edgelength,\qquad L(\vec \lambda,\vec \mu) \coloneqq \sum_{\node \in\setNode}\frac12 \sum_{\substack{\edge\sim \node\\ \edge =(\node,\node^\prime)}} \frac{(\vec \lambda_\node-\vec \lambda_{\node^\prime})\cdot (\vec \mu_\node-\vec \mu_{\node^\prime})}{\edgelength},
\end{align}
where the weighting with the edge length $\edgelength$ is chosen to be the same as for the mass and stiffness operators in a one-dimensional finite element implementation. 

The following theorem shows that the condensed problem \cref{EQ:hybrid_problem} is spectrally equivalent to the weighted graph Laplacians defined in \cref{EQ:graphlaplacian} in each component. This is done by explicitly computing the solutions to the one-dimensional Timoshenko beam equations, where we assume for simplicity that the material coefficients are edgewise constant. Under suitable assumptions on the discretized local solver~\cref{EQ:discrete_solver}, we expect that this result can be transferred to the condensed formulation of the HDG method in~\cref{EQ:hybrid_problem_disc}.

\begin{theorem}[Spectral equivalence to graph Laplacian]
	Assume that the maximal edge length is sufficiently small and that the material coefficients $C_{\vec n}$ and $C_{\vec m}$ are edgewise constant. Then, there holds for all $(\vec \lambda,\vec \phi) \in V_{\vec \lambda}\times V_{\vec \lambda}$ that
	\begin{align}
		\label{EQ:speceq}
		 L(\vec \lambda,\vec \lambda) + L(\vec \phi,\vec \phi) \lesssim A((\vec \lambda,\vec \phi),(\vec \lambda,\vec \phi)) \lesssim L(\vec \lambda,\vec \lambda) + L(\vec \phi,\vec \phi),
	\end{align}
	where the hidden constants depend only on the reciprocal of 
	\begin{align}
		\label{eq:friedrichs}
		\lambda_\mathrm{min} \coloneqq\min_{\vec \mu \in V_{\vec \lambda}\setminus \{0\}} \frac{L(\vec \mu,\vec \mu)}{M(\vec \mu,\vec \mu)},
	\end{align} 
which characterizes the smallest eigenvalue of the generalized eigenvalue problem with the bilinear forms $L$ and $M$ defined in~\cref{EQ:graphlaplacian} on the left- and right-hand sides, respectively. 
\end{theorem}

Note that $\lambda_\mathrm{min}$ can be uniformly bounded from below using results from classical graph theory; see~\cite{C05,CY95}. The bounds typically involve the constant of a $d$-dimensional isoperimetric inequality, which is a measure for the connectivity of the underlying graph. 
\begin{proof}
	The proof is done in two steps. First, for each beam, we explicitly solve the corresponding local Timoshenko beam equations, where the Dirichlet boundary data at the beam endpoints is given and the source terms are set to zero. Second, these explicit local Timoshenko solutions are used to derive the desired spectral equivalence result. 
	 
	\emph{Step 1:} We consider an arbitrary but fixed edge $\edge \in \setEdge$ and denote its associated local basis by~$\{\vec i_\edge,\vec j_\edge,\vec k_\edge\}$. As global basis we use the canonical basis of $\mathbb R^3$, denoted by $\{\hat{\vec i},\hat{\vec j},\hat{\vec k}\}$. To explicitly solve the local Timoshenko equations corresponding to~$\edge$ we write them in local coordinates similar to~\cref{rem:locform}. In the following, we denote local variables by hats and write~$\vec T_\edge$ for the change of basis matrix between local and global coordinates. Recall that the local coefficient matrices $\hat C_{\vec n}= \vec T_\edge^\top C_{\vec n}\vec T_\edge$ and $\hat C_{\vec m}= \vec T_\edge^\top C_{\vec m}\vec T_\edge$ are diagonal with positive diagonal entries, cf.~\cref{rem:locform}. The local Timoshenko equations are posed on the domain~$[0,\edgelength]\times \{0\}^{\textcolor{magenta}{2}}$, which we parametrize with the interval $[0,\edgelength]$. We omit the subscript $\edge$ to simplify the notation in the following. 
	 
	Let the following Dirichlet data for the displacement and rotation be given:
	 \begin{equation}
	 	\label{EQ:Dircond}
	 	\uhat(0) = \lamo,\qquad  \uhat(\edgelength) = \lamt, \qquad \rhat(0) = \phio,\qquad  \rhat(\edgelength) = \phit.
	 \end{equation}
 Similar to the definition of the local solver in~\cref{EQ:local_solver}, the local source terms are set to zero.  To derive the local Timoshenko solutions, we first note that $\partial_{\hat x} \nhat = 0$ and hence $\nhat = -\vec c$ for some $\vec c \in \mathbb R^3$. This implies for the displacement that $\partial_{\hat x} \uhat =- \ihat \times \rhat + \hat C_{\vec n}^{-1}\vec c$. To derive an explicit expression for the rotation $\rhat$ we note that $\partial_{\hat x}\mhat = -\ihat \times \nhat = \ihat \times \vec c$, which implies that $\mhat = (\ihat \times \vec c)\hat x-\vec d$ for some $\vec d \in \mathbb R^3$. The equation $\partial_{\hat x}\rhat = -\hat C_{\vec m}^{-1}\mhat$ can then be used to obtain an explicit expression for the rotation. Combining the above considerations yields the following expressions for the displacement and rotation:
\begin{align*}
	\uhat &= \tfrac16 (\ihat\times \hat C_{\vec m}^{-1}(\ihat \times \vec c))\hat x^3 - \tfrac12 (\ihat \times \hat C_{\vec m}^{-1}\vec d)\hat x^2 - (\ihat \times \phio)\hat x + \hat C_{\vec n}^{-1}\vec c \hat x + \lamo,\\
	\rhat &= -\tfrac12 \hat C_{\vec m}^{-1}(\ihat \times \vec c)\hat x^2 + \hat C_{\vec m}^{-1}\vec d \hat x +\phio.
\end{align*}
To determine the constants $\vec c$ and $\vec d$, we use the second and fourth conditions of \cref{EQ:Dircond}, which we have not yet incorporated in the above formula for $\uhat$ and $\rhat$. Below, we use the abbreviations $\lamdiff \coloneqq \lamt-\lamo$, $\phidiff \coloneqq \phit-\phio$, and $\phisum \coloneqq \phio+\phit$ to simplify the notation. For $\vec c$ we obtain, after some algebraic manipulations, the expression
\begin{align*}
	\vec c = \tfrac{1}{\edgelength}\hat C_{\vec n}\lamdiff + \tfrac12\hat C_{\vec n}(\ihat \times\phisum)  + \tfrac{\edgelength^2}{12}\hat C_{\vec n}(\ihat \times (\hat C_{\vec m}^{-1}(\ihat \times \vec c))),
\end{align*}
where $\vec c$ still appears also on the right-hand side. We explicitly evaluate the cross products to move the term involving $\vec c$ to the left-hand side. Introducing the matrix 
\begin{equation*}
	D\coloneqq \mathds{1}_3 +\tfrac{\edgelength^2}{12}\operatorname{diag} \big(0,(\hat C_{\vec n})_{22}(\hat C_{\vec m}^{-1})_{33},(\hat C_{\vec n})_{33}(\hat C_{\vec m}^{-1})_{22}\big),
\end{equation*}
where $(\cdot)_{ij}$ denote the $ij$-th entry of a matrix, this results in 
\begin{equation*}
	\vec c = \tfrac{1}{\edgelength}D^{-1}\hat C_{\vec n}\lamdiff + \tfrac12 D^{-1}\hat C_{\vec n}(\ihat \times \phisum).
\end{equation*}
This equation can then be used to derive the following explicit expression for $\vec d$:
\begin{equation*}
	\vec d = \tfrac{1}{\edgelength}\hat C_{\vec m}\phidiff + \tfrac12 (\ihat \times (D^{-1}\hat C_{\vec n}\lamdiff)) + \tfrac{\edgelength}4 (\ihat \times (D^{-1}\hat C_{\vec n}(\ihat \times \phisum))).
\end{equation*}
Inserting the above expressions for the constants $\vec c$ and $\vec d$ in the equations $\nhat = -\vec c$ and $\mhat = (\ihat \times \vec c)\hat x-\vec d$, we obtain for the forces and moments that
\begin{align*}
	\nhat &=-\tfrac{1}{\edgelength}D^{-1}\hat C_{\vec n}\lamdiff - \tfrac12 D^{-1}\hat C_{\vec n}(\ihat \times \phisum),\\
	\mhat &= \tfrac1{2\edgelength}(\ihat \times (D^{-1}\hat C_{\vec n}\lamdiff))(2\hat x-\edgelength) + \tfrac14(\ihat \times (D^{-1}\hat C_{\vec n}(\ihat \times \phisum)))(2\hat x-\edgelength) - \tfrac1{\edgelength}\hat C_{\vec m}\phidiff.
\end{align*}

To establish an expression for the bilinear form $A$ in terms of the Dirichlet data at the nodes, we will use its equivalent representation given in \cref{LEM:characterize_bil}. This representation includes integrals involving the forces and moments, which can be computed as follows:
\begin{align*}
	\int_{0}^{\edgelength} \hat C_{\vec n}^{-1}\nhat \cdot \nhat\,\mathrm{d}\hat x &= \tfrac1{\edgelength}|D^{-1}\hat C_{\vec n}^{1/2}\lamdiff|^2 + D^{-1}\hat C_{\vec n}^{1/2}\lamdiff \cdot D^{-1}\hat C_{\vec n}^{1/2}(\ihat \times \phisum) \\[-1ex]
	&\qquad + \tfrac{\edgelength}4 |D^{-1}\hat C_{\vec n}^{1/2}(\ihat \times \phisum)|^2,\\
	\int_{0}^{\edgelength} \hat C_{\vec m}^{-1}\mhat\cdot \mhat \,\mathrm{d}\hat x &= \tfrac{1}{\edgelength}|\hat C_{\vec m}^{1/2}\phidiff|^2 + \tfrac{\edgelength}{12}|\hat C_{\vec m}^{-1/2}(\ihat \times (D^{-1}\hat C_{\vec n}\lamdiff))|^2\\[-1ex]
	&\qquad + \tfrac{\edgelength}{48} |\hat C_{\vec m}^{-1/2}(\ihat \times (D^{-1}\hat C_{\vec n}(\ihat \times \phisum)))|^2\\
	&\qquad  + \tfrac{\edgelength^2}{12}\big(\hat C_{\vec m}^{-1/2}(\ihat \times (D^{-1}\hat C_{\vec n}\lamdiff))\big) \cdot \big(\hat C_{\vec m}^{-1/2}(\ihat \times (D^{-1}\hat C_{\vec n}(\ihat \times \phisum)))\big).
\end{align*}
To prove the lower bound in \cref{EQ:speceq}, we sum up the previously computed integrals and bound the result from below using the weighted Young inequality. This yields for any  $0<\epsilon<1$, which we will specify later, that
\begin{align*}
\int_{0}^{\edgelength} \hat C_{\vec n}^{-1}\nhat\cdot \nhat\, \mathrm{d}\hat x  &+ \int_{0}^{\edgelength}\hat C_{\vec m}^{-1}\mhat \cdot \mhat\, \mathrm{d}\hat x\\
 &\geq  \tfrac1{\edgelength}|D^{-1}\hat C_{\vec n}^{1/2}\lamdiff|^2 + \tfrac1{\edgelength}|\hat C_{\vec m}^{1/2}\phidiff|^2  +\tfrac{\edgelength}4 |D^{-1}\hat C_{\vec n}^{1/2}(\ihat \times \phisum)|^2\\
&\qquad + \big(D^{-1}\hat C_{\vec n}^{1/2}\lamdiff\big)\cdot \big(D^{-1}\hat C_{\vec n}^{1/2}(\ihat \times \phisum)\big)\\
&\geq \tfrac{1-\epsilon}{\edgelength}|D^{-1}\hat C_{\vec n}^{1/2}\lamdiff|^2 + \tfrac1{\edgelength}|\hat C_{\vec m}^{1/2}\phidiff|^2  -\tfrac{\edgelength(1-\epsilon)}{4\epsilon}|D^{-1}\hat C_{\vec n}^{1/2}(\ihat \times \phisum)|^2.
\end{align*}
Using the uniform bounds on the coefficients~\cref{eq:spectralbounds} and the assumption that the edge length~$\edgelength$ is sufficiently small, the latter estimate can be simplified to
\begin{equation}
	\label{EQ:estloc}
\int_{0}^{\edgelength} \hat C_{\vec n}^{-1}\nhat\cdot \nhat\, \mathrm{d}\hat x  + \int_{0}^{\edgelength}\hat C_{\vec m}^{-1}\mhat \cdot \mhat\, \mathrm{d}\hat x
	\gtrsim \tfrac{1-\epsilon}{\edgelength}|\lamdiff|^2 + \tfrac1{\edgelength}|\phidiff|^2-C\tfrac{\edgelength(1-\epsilon)}{4\epsilon}|\ihat \times \phisum|^2,
\end{equation}
where $C>0$ only depends on the bounds of the coefficients.

\emph{Step 2:} Estimate~\cref{EQ:estloc} was formulated in local coordinates, but it also holds in global coordinates, i.e., without hats. This is a direct consequence of the orthogonality of the change-of-basis matrix $\vec T_\edge$. Summing up the global estimates for all edges then yields  that
\begin{multline}
	\sum_{\edge \in \setEdge} \int_\edge C_{\vec n}^{-1}\vec n_\edge \cdot \vec n_\edge\ds  + \sum_{\edge \in \setEdge} \int_\edge C_{\vec m}^{-1}\vec m_\edge \cdot \vec m_\edge\ds \\
	\gtrsim  \sum_{\substack{\edge\in \setEdge\\ \edge =(\node,\node^\prime)}}\Big(\tfrac{1-\epsilon}{\edgelength}|\vec \lambda_\node-\vec \lambda_{\node^\prime}|^2 + \tfrac1{\edgelength}|\vec \phi_\node-\vec \phi_{\node^\prime}|^2 - C\tfrac{\edgelength(1-\epsilon)}{4\epsilon}|\vec i_\edge \times (\vec \phi_\node+\vec \phi_{\node^\prime})|^2\Big).\label{EQ:est}
\end{multline}
Rearranging the above sum into a sum over all nodes and another sum over all edges adjacent to the nodes, one observes that the first and second terms in \cref{EQ:est} can be written in terms of the bilinear form $L$ of the graph Laplacian defined in \cref{EQ:graphlaplacian}. Similarly, the last term in \cref{EQ:est} can be estimated using the mass-type operator $M$. This and using the definition of $\lambda_\mathrm{min}$ in \cref{eq:friedrichs} yields that 
\begin{align*}
		\cref{EQ:est} &\geq (1-\epsilon)L(\vec \lambda,\vec \lambda) + L(\vec \phi,\vec \phi) -C\tfrac{1-\epsilon}{\epsilon}M(\vec \phi,\vec \phi)\\
		& \geq (1-\epsilon)L(\vec \lambda,\vec \lambda) + (1-C\tfrac{1-\epsilon}{\epsilon}\lambda_\mathrm{min}^{-1})L(\vec \phi,\vec \phi).
\end{align*}
An investigation shows that for any $\epsilon \in \big(\tfrac{C}{\lambda_\mathrm{min}+C},1\big)$, the terms in the above estimate are positive, which proves the lower bound in \cref{EQ:speceq}.  The upper bound can be proved using similar arguments, and the proof is omitted for brevity.  This concludes the proof. 
\end{proof}
\section{Numerical experiments}\label{SEC:numexp}
%
In this section, we present numerical experiments that support the theoretical predictions of this paper. The numerical experiments have been implemented in the HyperHDG software package~\cite{HyperHDGgithub}, which is described in detail in \cite{RuppGK22}. There, numerical experiments demonstrating the convergence of the HDG method for diffusion-type problems on hypergraphs were presented. 
The code to reproduce the numerical experiments of this paper is available at \url{https://github.com/HyperHDG/}, and the data describing the fiber network model of paper used in the second numerical experiment can be found at \cite{data}.
\subsection*{Optimal order convergence of HDG method}
%

The first numerical example considers a toy problem to study the convergence properties of the proposed HDG method for Timoshenko beam networks. Specifically, we consider a network representing a two-dimensional unit cross embedded in the three-dimensional Euclidean space, i.e., $$\big([-1,1]\times \{0\} \cup \{0\}\times[-1,1]\big) \times \{0\}.$$ Before mesh refinement, the network consists of four edges and five vertices. Dirichlet boundary conditions are imposed at the four nodes located at the tips of the cross. The bounding domain is chosen to be $\Omega = [-1,1]^2\times \{0\}$. 
To construct suitable data and corresponding solutions, we use the method of manufactured solutions. More precisely, we consider homogeneous material coefficients, i.e., $C_{\vec n} = C_{\vec m} = \mathds 1$, and choose the force terms such that the problem admits the following displacements and rotations as its solutions:
\begin{equation*}
	\vec u(x,y,z) = \begin{pmatrix}
		0\\\cos(\pi y)\\
		\cos(\pi x)
	\end{pmatrix},\qquad 
\vec r(x,y,z) = \begin{pmatrix}
	0\\ \sin(\pi x) \\\sin(\pi y)
\end{pmatrix}.
\end{equation*}
\begin{figure}
	\begin{minipage}{.49\linewidth}
				\centering
		\includegraphics[width=.8\textwidth]{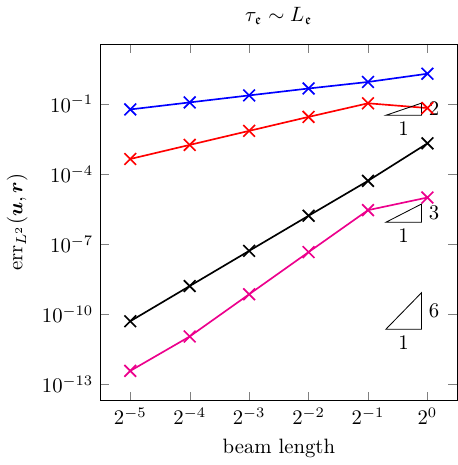}
	\end{minipage}
	\begin{minipage}{.49\linewidth}
				\centering
 \includegraphics[width=.8\textwidth]{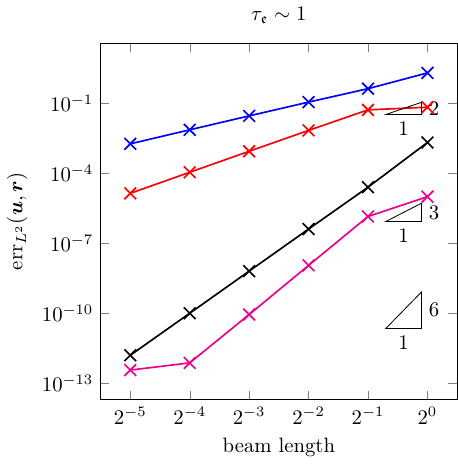}
	\end{minipage}\vspace{.5cm}
	\begin{minipage}{.49\linewidth}
						\centering
 \includegraphics[width=.8\textwidth]{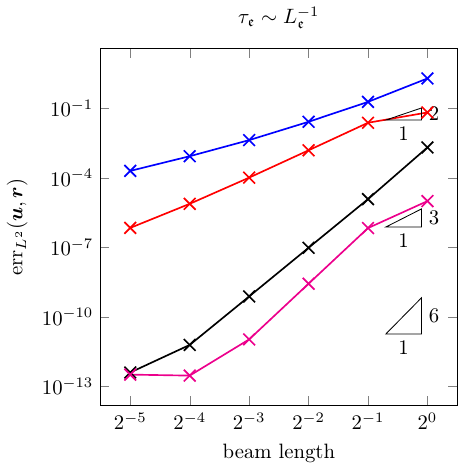}
 	\end{minipage}
 	\begin{minipage}{.49\linewidth}
 						\centering
  \includegraphics[width=.8\textwidth]{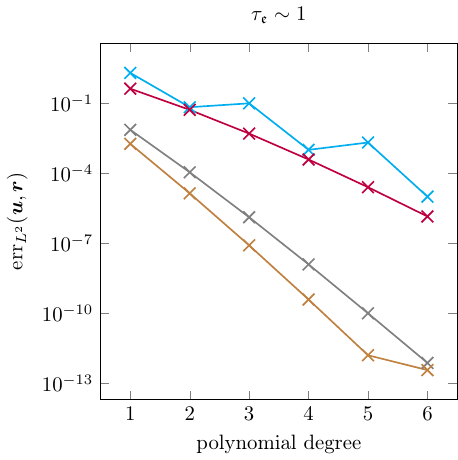}
   	\end{minipage}
 \caption{All four plots show the $L^2$-errors of the primal variables. The top left, top right, and bottom left plots show the errors for the polynomial degrees $1$ (blue), $2$ (red), $5$ (black), and $6$ (magenta) as a function of the mesh size, for different choices of the stabilization parameter $\tau$. The bottom right plot shows the errors for the beam lengths~$1$~(cyan), $2^{-1}$ (purple), $2^{-4}$ (gray), and $2^{-5}$ (brown) as a function of the polynomial degree.}\label{FIG:conv_plot}
\end{figure}

\cref{FIG:conv_plot} shows convergence plots for the $L^2$-error of the primal HDG variables, which we will subsequently abbreviate by $\operatorname{err}_{L^2}(\vec u,\vec r)$, cf.~\cref{EQ:L2primal}. In the top right and bottom left convergence plots, one observes optimal $(p+1)$-th order convergence for the classical choices of stabilization parameters $\tau_\edge \sim 1$ and $\tau_\edge \sim \edgelength^{-1}$. Note that for $\tau_\edge\sim \edgelength^{-1}$ one encounters a pre-asymptotic regime with even faster convergence. For the stabilization parameter $\tau_\edge\sim \edgelength$ one observes a suboptimal $p$-th order convergence, cf.~\cref{FIG:conv_plot} (top left). These convergence results are consistent with the theoretical predictions of \cref{COR:errorest}. Note that the convergence rates we obtain are also consistent with those found in the literature. Convergence tests for the dual variables also confirm the rates predicted by \cref{COR:errorest}. For the sake of brevity, however, the corresponding plots are not shown here. In addition to the convergence under mesh refinement, we also want to numerically investigate the convergence as the polynomial degree is increased. The corresponding plot can be found at the bottom right in \cref{FIG:conv_plot}. One observes an exponential convergence as the polynomial degree is increased, with the rate of decay depending on the fixed mesh size. A similar exponential convergence behavior can also be observed for the dual variables.

\subsection*{Elastic deformation of paper}
%

\begin{figure}
\begin{minipage}{.49\linewidth}
	\centering
\includegraphics[width=.95\textwidth]{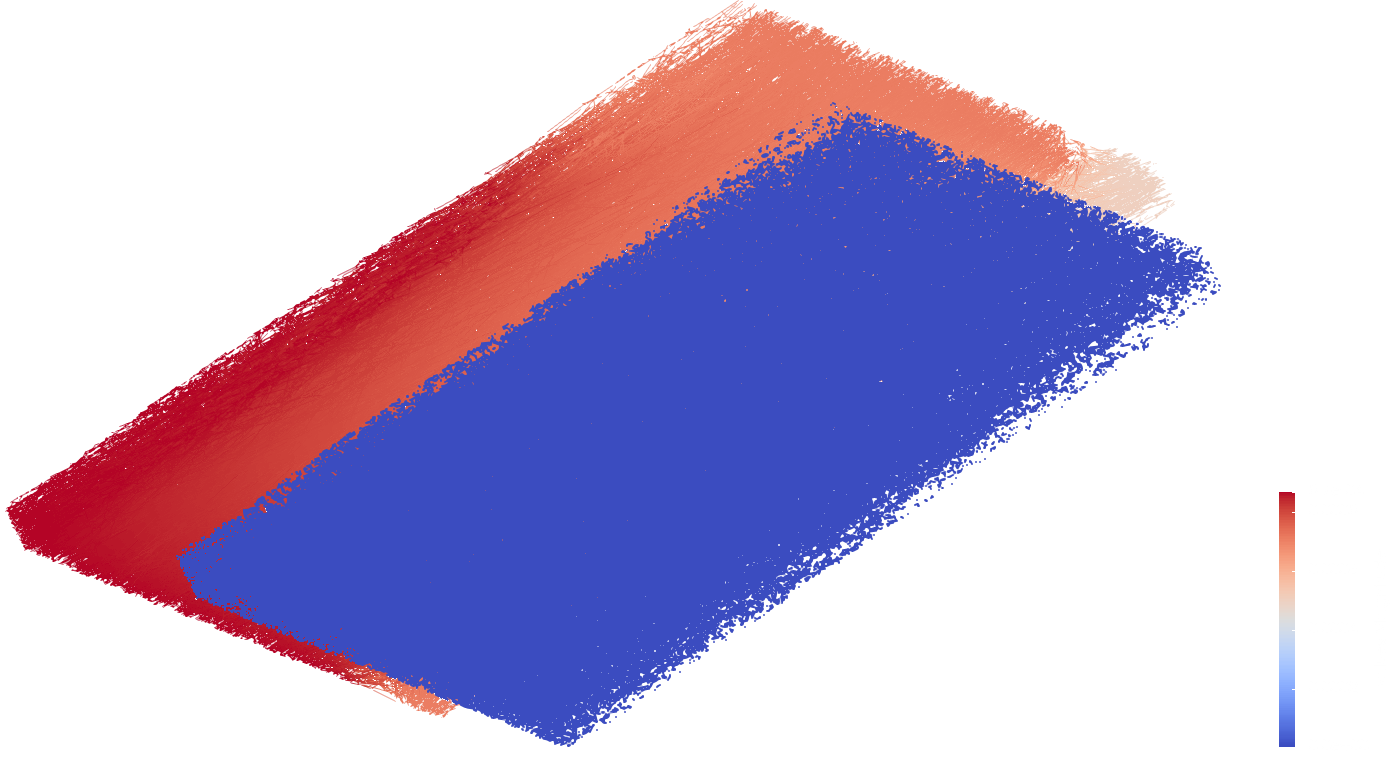}
\end{minipage}\begin{minipage}{.49\linewidth}
\centering
\includegraphics[width=.8\textwidth]{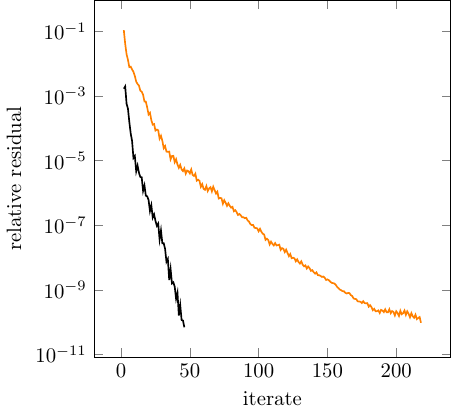}
\end{minipage}
 \caption{On the left an illustration of the network before (blue) and after (red) deformation is shown. On the right is  a convergence plot of the preconditioned conjugate gradient method. There the relative residual is plotted as a function of the iteration number for constant material parameters (black) and realistic material parameters (orange).}\label{FIG:real_paper}
\end{figure}

The purpose of the second numerical experiment is to demonstrate the applicability of the proposed HDG method to a realistic example. We consider the elastic deformation of about 2 mm x 4 mm piece of paper, where our collaborators at the Fraunhofer-Chalmers Centre (FCC) provided the corresponding spatial network and material parameters. The spatial network consists of about 615K edges and 424K nodes. We consider the stretching of the paper caused by inhomogeneous Dirichlet boundary conditions at nodes at the lateral boundary. We use a HDG discretization with polynomial degree 5 for all edges. The resulting linear system of equations is then solved by the preconditioned conjugate gradient method with the domain decomposition preconditioner as described in \cref{SEC:precond}. For the preconditioner, we use an artificial coarse mesh consisting of 64 elements, where only one element is used in the $z$-direction. The local subproblems arising within the preconditioner are solved with a conjugate gradient method without preconditioning using a termination criterion of a relative residual of $10^{-3}$.  This is practically feasible since we use a mild termination criterion. 
The computations were performed in parallel on a cluster with 64 cores. For the case of constant material parameters, our solver needed 46 iterations and about 9 hours to solve the problem up to a relative residual of $10^{-10}$. For the realistic parameters, it took 218 iterations and about 33 hours to reach the same accuracy. The resulting deformed piece of paper and the convergence plot of the preconditioned conjugate gradient method are shown in \cref{FIG:real_paper}.
 
\section{Acknowledgments}

We would like to thank Morgan Görtz (Fraunhofer-Chalmers Center) for many fruitful discussions on the application of fiber-based materials, for providing some of the figures in this paper, and for generating the network data used in the second numerical experiment. 
%
\bibliographystyle{ARalpha}
\bibliography{hdg_network_references}
%
\end{document}